\documentclass[final]{elsarticle}

\usepackage{amsmath,amssymb}
\usepackage{graphicx}
\usepackage{footnote}
\usepackage{algorithmic}
\graphicspath{{./images/}}
\usepackage{natbib}
\usepackage{xcolor}

\title{Grazing bifurcations and transitions between periodic states of the PP04 model for the glacial cycle.}
\author[bath]{Chris~Budd*}
\ead{mascjb@bath.ac.uk}
\author[biust]{Kgomotso S. Morupisi}
\ead{morupisik@biust.ac.bw}
\cortext[cor1]{Corresponding author}
\address[bath]{Department of Mathematical Sciences, University of Bath, BA2 7AY, UK}
\address[biust]{Department of Mathematics and Statistical Sciences, Botswana International University of Science and Technology,Private Bag 16, Palapye, Botswana}

\begin{document}

\begin{abstract}
    We look at the periodic behaviour of the Earth's glacial cycles and the transitions between different periodic states when either external parameters (such as $\omega$) or internal parameters (such as $d$) are varied. We model this using the PP04 model of climate change. This is a forced discontinuous Filippov (non-smooth) dynamical system. When periodically forced this has co-existing periodic orbits. We find that the transitions in this system are mainly due to grazing events, leading to grazing bifurcations. An analysis of the grazing bifurcations is given and the impact of these on the domains of attraction and regions of existence of the periodic orbits is determined under various changes in the parameters of the system. Grazing transitions arise for general variations in the parameters (both internal and external) of the PP04 model. We find that the grazing transitions between the period orbits resemble those of the Mid-Pleistocene-Transition. 
\end{abstract}

\maketitle

\section{Introduction}
\noindent One of the most interesting problems concerning the behaviour of the Earth's climate, is an understanding of the glacial cycles  over the last few million years. Measurements of the Earth's temperature over this period can be obtained by studying the concentration of the Oxygen-18 isotope in ice cores taken from Antarctica \cite{petit1999climate}. These shows that from around 2.5 Myr ago, the temperature had moderate sized oscillations with a period of around 40 kyr. However, at about 1 Myr ago, these changed into much larger amplitude oscillations with a period close to 100 kyr.

\vspace{0.1in}
 
 \noindent 
The  sudden change of  glacial oscillations  from low amplitude symmetric cycles  with a period of  about  $40$ kyr to larger amplitude asymmetric cycles with a  period of $100$ kyr, is known as Mid-Pleistocene Transition (MPT).
The reason for this behaviour is unclear, as during this transition there was no significant  change in the (Milankovitch) astronomical forcing  \cite{ashwin2015middle}. Several hypotheses have been proposed to explain the mechanisms behind the glacial cycles and the MPT. These are usually a combination of models for the interaction between the astronomical forcing and the internal dynamics of the Earth's climate \cite{ditlevsen2018complex}. A number of dynamical system models have been proposed for this, see for example \cite{saltzman1991first}, \cite{saltzman1990first}, \cite{jones2005global}, \cite{docquier2019impact}, \cite{widiasih2013dynamics}. Many of these which describe the Earth's climate using  smooth dynamical system models and link the behaviour to that of an excitable system. It is then hypothesised that the periodic glacial cycles arise via Hopf or cyclic fold bifurcations, and that the MPT arises from slow changes in the parameters of these systems (due for example to tectonic drift). An alternative hypothesis is that the glacial cycles are best described by relaxation oscillators, with rapid transitions between states \cite{crucifix2012oscillators}, \cite{ashwin2015middle}, \cite{paillard1998timing}. One motivation for this idea is that 
rapid changes in  global Carbon Dioxide, as observed from  proxy data, could be the cause of the rapid transition between the inter-glacial and the glacial states. That is, the growth of the ice sheets is terminated by a sudden growth of Carbon Dioxide in the atmosphere leading to an inter glacial state. The  gradual decline of this Carbon Dioxide to a steady state (typically through absorption into the oceans) reduces the greenhouse effect, which promotes the build up of the Northern Hemisphere ice sheets and a further period of glaciation \cite{willeit2019mid}. Based on this observation Paillard  \cite{paillard2004antarctic} proposed that the glacial cycle transition occurs when  a certain {\em threshold} value of the bottom salty water formation (Ocean stratification) is reached, promoting a sudden ventilation of Carbon Dioxide into the atmosphere.  It was shown in \cite{crucifix2012oscillators}, \cite{morupisi2020analysis} that when coupled to periodic astronomical forcing  this model (referred to as PP04) leads to periodic glacial cycles (which are only slightly perturbed when the astronomical forcing is quasi-periodic).  Mitsui \textit {et al} \cite{mitsui2015bifurcations} asserts that parameter sensitivity in some climate  models can be associated with  thresholds with some discontinuity  surface, such as the sudden ventilation described above. It is studying this phenomenon, using the framework of bifurcations in non-smooth dynamical systems, which is the motivation behind this current paper.

\vspace{0.1in}

\noindent The PP04 climate model is an example of a  forced piece-wise smooth dynamical system of Filippov type, characterised by a phase space  divided into two regions (glacial and interglacial state), by a discontinuity boundary $\Sigma$. In each region, the evolution of the trajectories are defined by different smooth dynamical systems, but there is a discontinuous change between these across $\Sigma$. In the PP04 model this discontinuity corresponds to the sudden ventilation of the Carbon Dioxide from the Southern Ocean. Such Filippov systems can have fixed points and periodic orbits. Significant changes in the dynamics can arise as parameters are varied, when the fixed points intersect $\Sigma$  (border collision bifurcations) or the periodic orbits have non-transversal intersections with $\Sigma$ (grazing bifurcations). Both of these are examples of the thresholding phenomena identified in \cite{mitsui2015bifurcations} and the different types of transitions at them are described in \cite{bernardo2008piecewise}. Characteristics of such transitions are sudden changes from fixed points to periodic orbits (at border collision bifurcations), and equally sudden changes from one form of periodic orbit to another (at a grazing bifurcation). We will show in this paper that grazing bifurcations play a very important role in explaining the observed dynamics of the PP04 model giving a {\em generic} mode of transition from one state to another as the  parameters in the model are varied. Note that these parameters can either be 'external' such as the frequency and amplitude of the insolation, or 'internal' such the salty water bottom formation efficiency controlled by the parameter $d$. Furthermore, we see behaviours  both in the evolution of transient states and as model parameters vary which resemble those seen at the MPT. An example of this behaviour is illustrated in Figure \ref{drange} in Section 3.3, where we see a change from a low amplitude periodic solution, to a larger amplitude periodic solution of higher period. We have no direct evidence that changes in these parameters led to the behaviour observed at the MPT. However, the generic nature of the grazing bifurcation, and the observation that under general parameter variation we see the transitions from low amplitude, low period periodic solutions to high amplitude, high period periodic solutions, suggests that the MPT may have arisen from a grazing bifurcation under a parameter change. 
 
\vspace{0.1in}

\noindent The structure of the remainder of this paper is as follows. In Section 2 we briefly review the derivation of the PP04 model. In Section 3 we review some of the dynamics of the periodically forced PP04 model and give examples of the transitions seen in periodic orbits at a grazing bifurcation (together with associated bifurcation diagrams). In Section 4 we derive the form of the Poincare map close to grazing. In Sections 5 and 6 we derive the analytical form of the grazing surface ${\cal G}$ (the set of initial values which lead to orbits which have non-transversal intersections with $\Sigma$), and then compute it for certain values of the parameters. In Section 7 we show how the analytical form of ${\cal G}$ is strongly reflected in the observed forms of the domains of attraction, with orbits starting close to the boundaries of the domains of attraction having changes from transient behaviour resembling one periodic orbit to limiting behaviour resembling another, occurring at grazing events. In Section 8 we show how the regions of existence of the various periodic orbits are shaped by the curves (in parameter space)  of the grazing bifurcations, and how this depends upon the various parameters in the PP04 model. In Section 9 we consider the impact of quasi-periodic forcing on the PP04 model, and show that this is also very much influenced by grazing events. Finally in Section 10 we draw some conclusions.

\section{The PP04 model for the ice ages and its basic dynamics}

\subsection{The PP04 model}

\noindent In this section we summarise the results presented in \cite{morupisi2020analysis}, and in \cite{paillard2004antarctic}. The PP04 model for the ice ages describes the system state through  three piece-wise linear ordinary differential equations in the total ice volume $V$, the Antarctic ice $A$ and the CO2 level $C$ forced by (quasi-periodic) Solar insolation. This model describes the Earth's climate as having two states, namely a cold glacial state and a warm inter-glacial state. The model assumes that there is a sudden transition between these two states when a large amount of Carbon Dioxide $C$ is released from the Southern Ocean into the atmosphere, leading to a general warming. The cause of this release is a change in the stratification of the ocean, which is triggered by an interplay between the global ice volume $V$ and the Antarctic ice cover $A$. The PP04 model is described as follows:
    
\begin{eqnarray}
%\begin{align}
%\begin{split}
  \frac{dV}{dt}& = &\frac{1}{\tau_V}(-x C - y I_{65}(t) + z - V),\\
  \frac{dA}{dt}& =& \frac{1}{\tau_A}(V-A),\\
  \frac{dC}{dt}& =& \frac{1}{\tau_C}( \alpha I_{65}(t) - \beta V + \gamma H(-F) + \delta -C).
%\end{split}
%\end{align}
\label{cjan1}
\end{eqnarray}
Here
\begin{equation}
    F = aV-bA - cI_{60}(t)+d
\label{cjan2}
\end{equation}
is called the salty bottom water formation efficiency (oceanic stratification) parameter. Of particular significance in this model is the discontinuity introduced by the Heaviside function $H(-F)$. In Crucifix \cite{crucifix2012oscillators} a smoothing of this is considered in which we replace $H$ by the smoother (but still significantly nonlinear) sigmoid function
\begin{equation}
 H_{\eta}(x) = \frac{1}{2} \left( 1 + \tanh(\eta x) \right) = \frac{1}{1 + e^{-2\eta x}}
 \label{H1}
 \end{equation}
For the numerical calculations presented in this paper we will usually use this function with $\eta = 1500$, which gives a close approximation to the discontinuous case. We will also use the default \cite{paillard2004antarctic} values of $x=1.3, y =0.5, z = 0.8, \alpha = 0.15, \beta =0.5,\gamma = 0.7,  \delta =0.7 , a= 0.3,b=0.7,d =0.27.$ We will also consider the effect of varying some of these parameters, in particular the frequency and amplitude of the forcing, and the parameter $d$ in the definition of $F$. 

\vspace{0.1in}

\noindent The forcing function,  $I_{65}(t)$, is the quasi-periodic Solar insolation forcing at $65^{\circ}$ North. For the majority of this paper we will simplify this to  be the periodic function
$$I_{65} = \mu \sin(\omega t).$$
This is a reasonable approximation, given the actual form of $I_{65}$ discussed in \cite{morupisi2020analysis}. It also allows us to make useful analytic predictions. In Section 9 we will consider the case of quasi-periodic forcing and show that an understanding of grazing behaviour allows us to also make predictions of the dynamics in this case.

\vspace{0.1in}

\noindent The piece-wise linear ODE given by the PP04 model admits global in time solution trajectories ${\mathbf X}(t)$ which typically converge towards synchronised period-solutions, modelling glacial cycles. It is the generic form of the transitions between these periodic states as parameters in the PP04 model (such as $\mu, \omega$ and $d$ are varied) which are of interest to us in this paper.

\subsection{Formulation of the PP04 model as a discontinuous Filippov system}

\noindent For the development of the theory in this paper we formulate the PP04 model (\ref{cjan1}) as a parametrised discontinuous Filippov System. In particular we have a transition surface $\Sigma$,  
across which the dynamics changes. The interaction between ${\mathbf X}(t)$ and $\Sigma$ is at the heart of the theory developed in this paper. In \cite{morupisi2020analysis} (see also \cite{paillard2004antarctic}), it is argued that we can neglect the $I_{60}(t)$ contribution to (\ref{cjan2}). We then have this surface is given by the zero set of the linear function
\begin{equation}
{
F({\mathbf X}) \equiv {\mathbf c}\cdot{\mathbf X} + d.
}
\label{c1}
\end{equation}

We define the following two states corresponding to the glacial and interglacial periods:
\begin{equation}
{
S^{+} = \{ {\mathbf X} : F({\mathbf X})> 0 \} , \quad S^{-} = \{ {\mathbf X} : {F(\mathbf X}) < 0. \}.
}
\label{c2}
\end{equation}
The PP04 model can then be expressed (see \cite{morupisi2020analysis}) as a discontinuous Filippov system in the form:
\begin{equation}
{
\dot{{\mathbf X}} = L {\mathbf X} + {\mathbf b}^{\pm} + I_{65}(t) \; {\mathbf e}.
}
\label{c3}
\end{equation}
Here the vectors ${\mathbf b}^{\pm}$ depend on the system state, and the linear operator $L$ and the vector ${\mathbf e}$ do not depend upon the system state.

\vspace{0.1in}

\noindent The smoothed PP04 model then takes the form
\begin{equation}
{
\dot{{\mathbf X}} = L {\mathbf X} + {\mathbf b}^{-} +
H_{\eta}(F) ({\mathbf b}^+ - {\mathbf b}^-) + I_{65}(t) \; {\mathbf e}.
}
\label{c3aa}
\end{equation}

\subsection{Trajectories}

\noindent The PP04 model is a piece-wise linear system, which has linear behaviour between intersections of the trajectories with $\Sigma$.
 The linearity allows an exact solution of the dynamics between these events. In particular, if ${\mathbf X}(t_k) = {\mathbf X}_k$ then we have

\begin{equation}
    {\mathbf X}(t) = e^{L(t-t_k)} \left({\mathbf X}_k + L^{-1} {\mathbf b}^{\pm} -{\mathbf c}(t_k) \right)  - L^{-1} {\mathbf b}^{\pm} + {\mathbf p}(t).
    \label{H3}
\end{equation}

\noindent Here ${\mathbf p}(t)$ is a {\em particular solution} of the system
$$\dot{{\mathbf p}} = L {\mathbf p}  +  I_{65}(t) \; {\mathbf e}.$$

\vspace{0.1in}

\noindent {\bf Lemma 2.1} {\em If $I_{65}(t) = \sum_n a_n e^{i \omega_n t}$ then
\begin{equation}
    {\mathbf p}(t) =  \sum_n a_n (i \omega_n - L)^{-1} e^{i \omega_n t} {\mathbf e}.
    \label{H2}
\end{equation}
}

\vspace{0.1in}

\noindent {\em Proof} This follows from a direct computation. Note that by taking real and imaginary parts we can readily compute the solutions to forcing functions
$\sin(\omega t)$ and $\cos(\omega t).$ \qed

\vspace{0.2in}

\noindent A solution of (\ref{cjan1}) then comprises a series of solutions of the form (\ref{H3}) with impacts with the surface $\Sigma$ at the times $t_k$, at which points ${\mathbf X}$ is continuous,  $\dot{\mathbf X}$ has a discontinuity, and the expression in (\ref{H2}) alternates between taking ${\mathbf b}^+$ and ${\mathbf b}^-$ in its argument. 

\vspace{0.1in}

\noindent A useful measure of the trajectory is the value of the function $F(t)$. In general $F$ is a smooth function, but it loses smoothness in its second derivative at the intersection points. The following result is proved in \cite{morupisi2020analysis}

\vspace{0.2in}

\noindent {\bf Lemma 2.2} {\em At an intersection point $t_i$ in which the trajectory passes from $S^+$ to $S^-$ (or vice-versa) there is a constant $r$ for which

\begin{equation}
[F]^+_- = [\dot{F}]_-^+ = 0, \quad [\ddot{F}]^+_- = r > 0.
\label{Fdisc}
\end{equation}
}

\vspace{0.2in}

\subsection{Periodic orbits and their ranges of existence}

\noindent It is shown in \cite{morupisi2020analysis} that if we take $I_{65}(t)$ to be periodic,
so that (for example)
$$I_{65}(t) = \mu \sin(\omega t),$$
then for certain values of $\mu$ and $\omega$ the PP04 model admits synchronised periodic solutions. These solutions are described in \cite{morupisi2020analysis} as $(m,n)$ orbits if they have $m$ glacial cycles, and period $T = 2 \pi n/\omega.$ 

\vspace{0.1in}

\noindent If $\mu = 0$ (the system is unforced) then for a range of value of the parameters, the system has an attracting {\em periodic orbit} of period $P \equiv 2\pi/\omega^*$.
For the parameter values in \cite{paillard2004antarctic} we have $P = 143 kyr, \omega^* = 0.0439.. .$ 

\vspace{0.1in}

\noindent It is shown further in \cite{morupisi2020analysis} that if $\mu \ll 1$ then the $(m,n)$ orbits exist in 'tongues' of parameter values $(\omega,\mu)$ so that as
$\mu \to 0$ we have that 
\begin{equation}
    \omega \to \omega_{m,n} = \frac{n \omega^*}{m} \quad \mbox{with} \quad |\omega - \omega_{m,n}| = {\cal O}\left(\mu^m\right). 
    \label{cjan3}
    \end{equation}
For these small values of $\mu$ the boundaries of the tongues correspond to saddle-node or period-doubling bifurcations.

\vspace{0.1in}

\begin{figure}[htb!]
\centering
\includegraphics[scale=0.125]{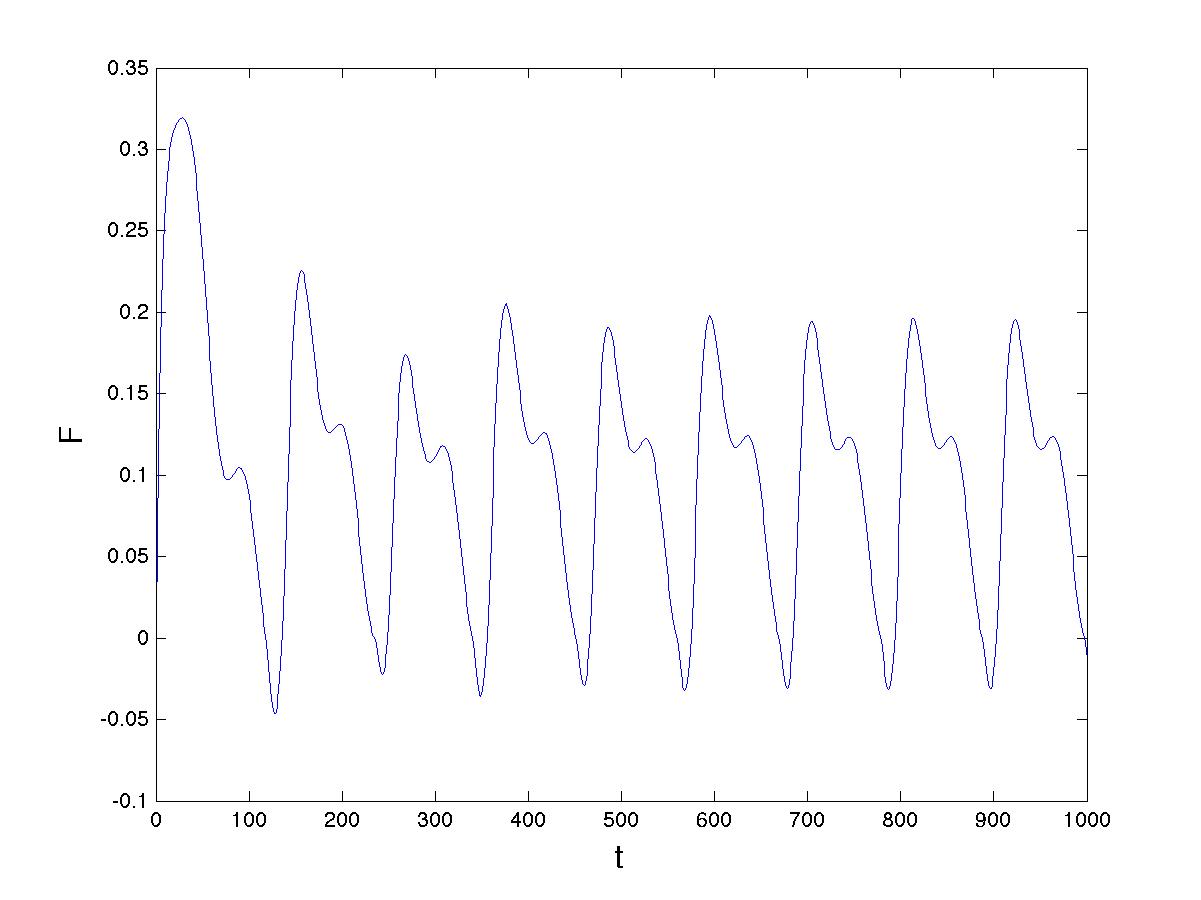}
\includegraphics[scale=0.125]{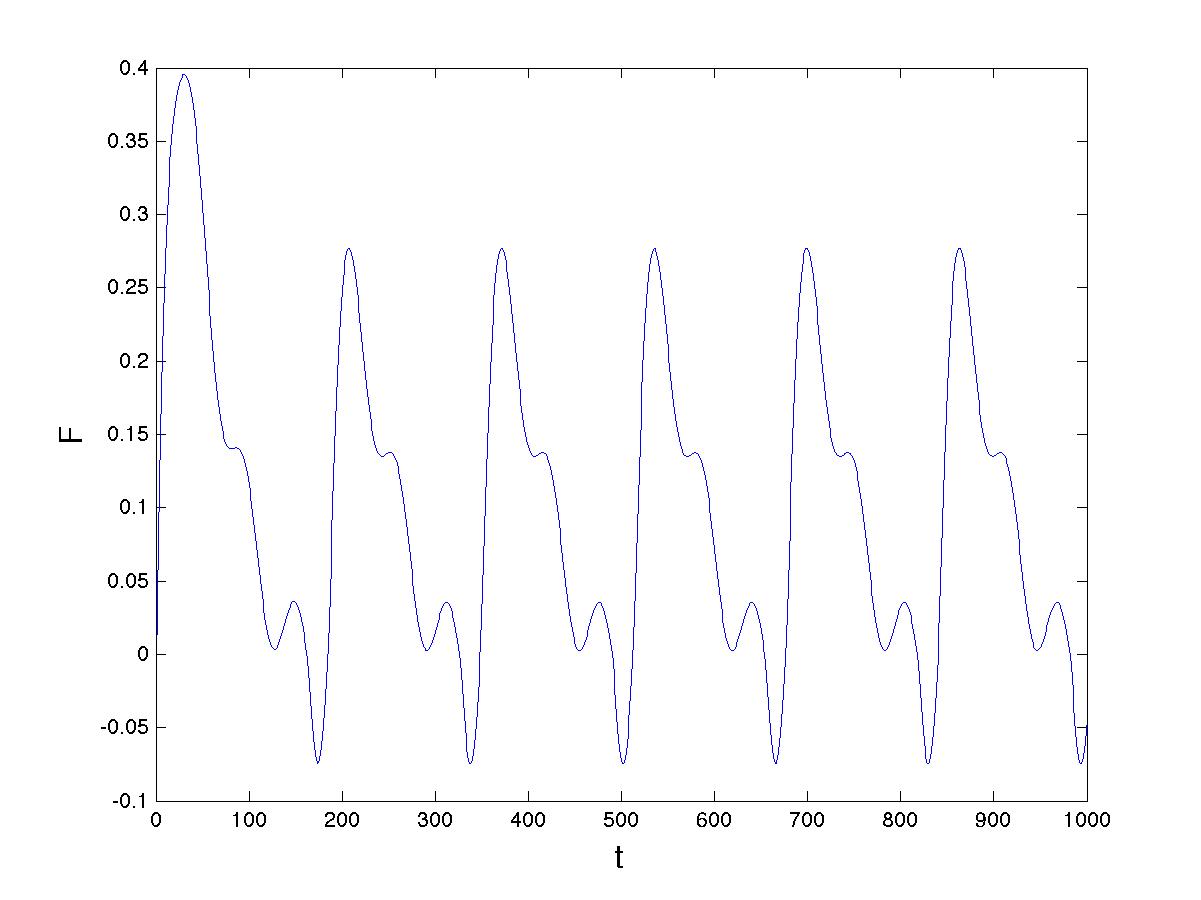}
\caption{Co-existing $(1,2)$ and $(1,3)$ periodic orbits at $\mu = 0.3$ $\omega = 0.115$, with  $\eta = 1500$.} 
\label{fig:ped2}
\end{figure}

\noindent For larger values of $\mu$ (for example $\mu > 0.2$) we see an overlapping of the tongues with the simultaneous existence of different stable periodic orbits. \noindent When, for example, $\mu = 0.3$ and $\omega = 0.115.$ we see co-existing $(1,2)$ and $(1,3)$ periodic solutions.  In Figure \ref{fig:ped2} we plot these periodic orbits, and an initial transient, showing $F$ as a function of $t$. This allows the form of the orbit and its intersection with $\Sigma$ to be seen most easily. Each periodic orbit has a {\em glacial region} with $F > 0$ and an 
{\em inter-glacial region} with $F < 0$. The transitions between these regions occur when $F = 0$. We can see a number of additional oscillations within the glacial region (of size proportional to $\mu$). 

\vspace{0.1in}

\begin{center}
\begin{figure}[htb!]
\centering
\includegraphics[scale=0.28]{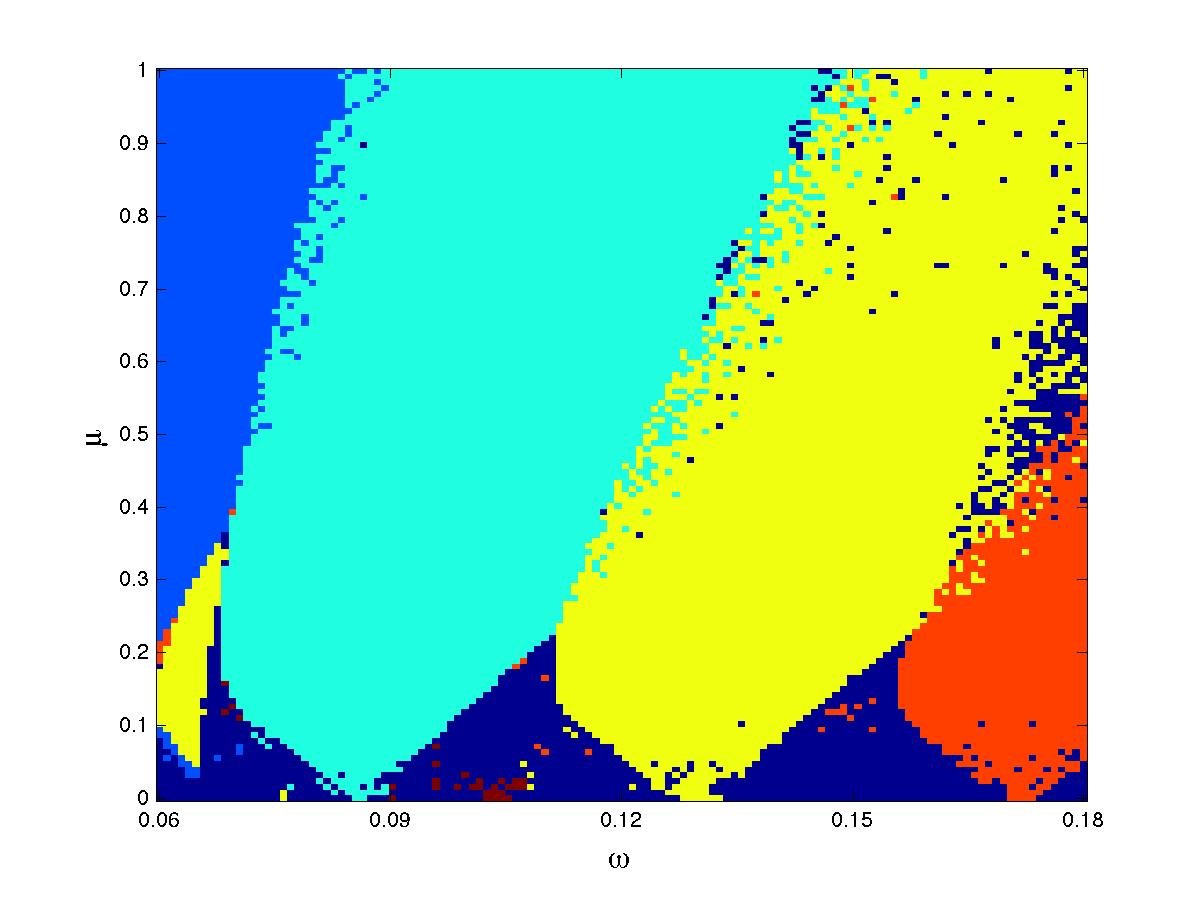}
\caption{The existence regions for the $(1,1)$ (light blue),$(1,2)$ (cyan),$(1,3)$ (yellow),$(1,4)$ (red) and $(2,5)$ (maroon), periodic orbits together with quasi-periodic behaviour (dark blue). For $\mu > 0.2$ the existence regions overlap.} 
\label{fig:fig1}
\end{figure}
\end{center}

\noindent A calculation of the existence regions of the $(m,n)$ orbits for the 'standard' set of parameters given in \cite{paillard2004antarctic} is given in Figure \ref{fig:fig1} in which we illustrate the (co)-existence of $(1,1), (1,2),(1,3),(1,4)$ and $(2,5)$ orbits as well as quasi-periodic behaviour. This figure was computed by taking random initial data for the system (\ref{c3aa}),  with $\eta = 1500.$ The resulting $\omega-$limit set is then computed using the Matlab code {\tt ode15s} to solve the smoothed system (\ref{c3aa}) for a long period, and the nature of the periodic orbit found from inspecting the peaks of an {\tt fft} of the solution. 

\vspace{0.1in}

\noindent We see clearly in this figure that for small $\mu$ the existence regions have the 'Arnold Tongue' structure predicted by (\ref{cjan3}) with the regions bounded by smooth bifurcations. This is typical of smooth dynamical systems and follows for the PP04 model as the trajectories for small $\mu$ intersect $\Sigma$ transversely. However this structure breaks down for larger values of $\mu$ when we can have non-transversal 
{\em grazing} intersections with $\Sigma$. These are then associated with {\em grazing bifurcations} and lead to dramatic change in the solution behaviour. We speculate in this paper both that such {\em non-smooth bifurcations} may be linked to the MPT as certain parameters in the PP04 model (such as $d$) vary,  and also that they account for some of the sensitivity observed in climate dynamics. Such behaviour was noted in \cite{paillard2004antarctic} (see also  \cite{mitsui2015bifurcations}) where it was described as 'thresholding'. We now study it in detail.

%\begin{center}
%\begin{figure}[htb!]
%\includegraphics[scale=0.35]{mu01.jpg}
%\caption{Solution tongues showing a $(1,2)$ $(1,3)$ and %orbit together with quasi-periodic behaviour. Note the %transitions close to $mu=0.3$ and $\omega = 0.115.$} 
%\label{fig:fig2}
%\end{figure}
%\end{center}

\vspace{0.1in}

\section{Transitions as parameters vary}

\subsection{Sudden grazing transitions}

\noindent Our primary interest in this paper lies in the possible transitions as parameters vary, between the different forms of periodic orbit as illustrated, for example,  in Figure \ref{fig:fig1}. In particular we will consider changes in both $\omega$ and $d$. 

\vspace{0.1in}

\noindent {\em Varying $\omega$.}

\vspace{0.1in}

\noindent As a first 
calculation we compute a Mont\'e-Carlo bifurcation diagram of the existence regions for fixed $\mu$ and $d$ and varying $\omega$. To do this we initially keep $\mu = 0.3, d= 0.27$ constant and slowly increase $\omega$ from $\omega = 0.01$ to $\omega = 0.15$. For each value of $\omega$ we take 10 random initial states and compute the resulting trajectories. We then plot the $\omega-$limit set of the values of the function $F$.  The resulting bifurcation diagram is shown in Figure \ref{fig:fig3a}
\begin{center}
\begin{figure}[htb!]
\includegraphics[scale=0.25]{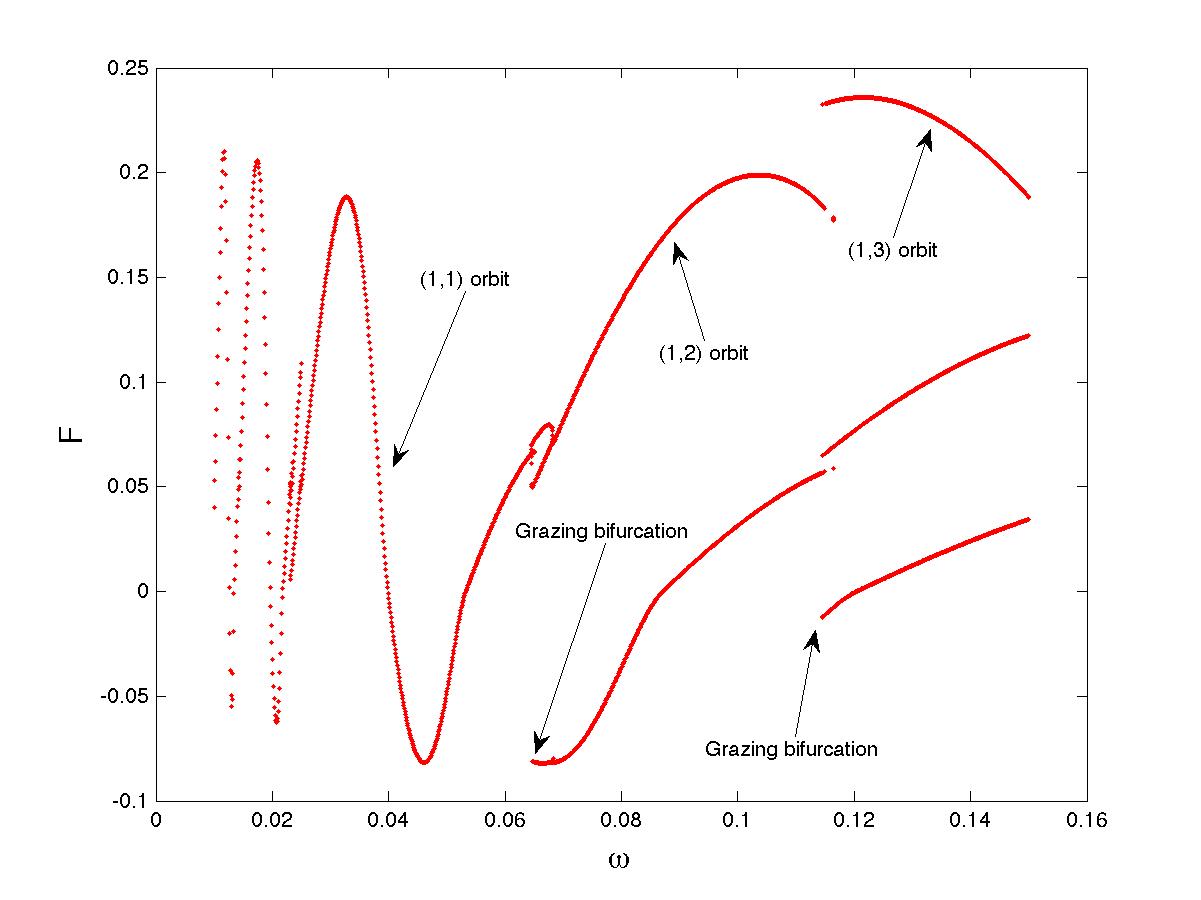}
\caption{Mont\'e-Carlo bifurcation diagram when $\mu=0.3$ and $0.01 < \omega < 0.15.$ In this we see abrupt transitions between the $(1,1)$, $(1,2)$ and $(1,3)$ orbits at grazing bifurcations.}
\label{fig:fig3a}
\end{figure}
\end{center}
In this figure we see a $(1,3)$ orbit for the larger values of $\omega$. As $\omega$ is decreased, this then vanishes abruptly at $\omega \approx 0.114$ where it is replaced by a $(1,2)$ periodic orbit. This then vanishes as $\omega$ is decreased at $\omega \approx 0.065$ and is replaced by a $(1,1)$ orbit. The $(1,1)$ orbit then persists as $\omega$ decreases further, but exhibits complex behaviour for the smaller values of $\omega$ (which for the climate problem are not physically interesting values.) The abrupt changes at $\omega \approx 0.065$ and $\omega \approx 0.115$ are both examples of {\em grazing bifurcations}.
To illustrate this in more detail we now show the bifurcation diagram for $\omega$ close to $\omega = 0.115$ in Figure \ref{fig:fig3}.
\begin{figure}[htb!]
\centering
\includegraphics[scale=0.15]{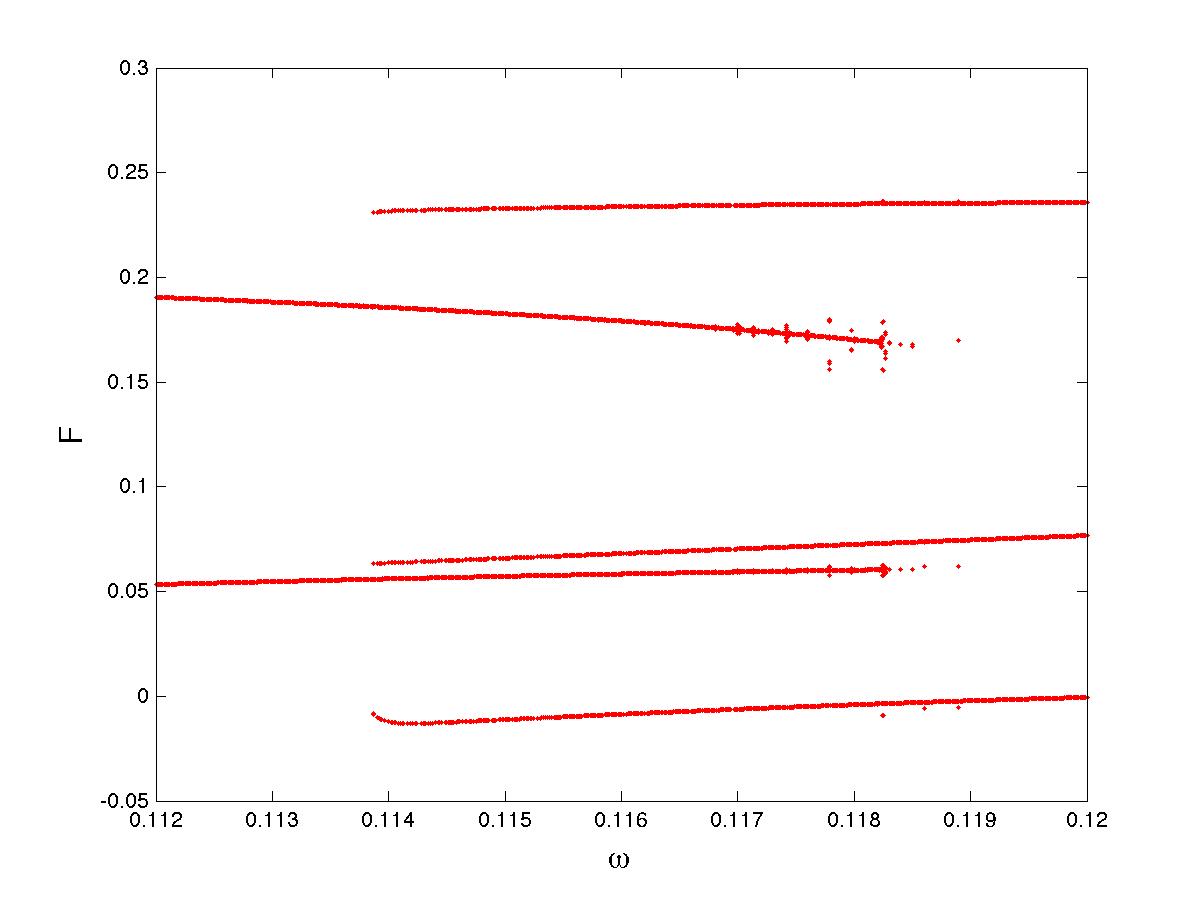}
\caption{Mont\'e-Carlo bifurcation diagram when $\mu=0.3$ and $\omega \approx 0.115$. In this we see a $(1,2)$ orbit which vanishes at a period-doubling bifurcation at $\omega \approx 0.118$ and a  $(1,3)$ orbit which vanishes at a {\em grazing} bifurcation at $\omega \approx 0.114$. }
\label{fig:fig3}
\end{figure}
Here we see a $(1,3)$ periodic orbit for the larger values of $\omega$ and a $(1,2)$ periodic orbit for the smaller values. As $\omega$ {\em increases} the $(1,2)$ orbit exhibits a smooth super-critical {\em period-doubling} bifurcation at $\omega \approx 0.118$.
Following the bifurcation there is a $(2,4)$ orbit close to the $(1,2)$ orbit, followed by what appears to be a period-doubling cascade leading to a small region of chaotic behaviour. Such a bifurcation structure is commonly observed in many smooth dynamical systems. The {\em right hand boundary} of the region of existence of the $(1,n)$ periodic orbits for the larger values of $\mu$ seems, from computations, always to be determined by such a period-doubling bifurcation. (For example, we can see in Figure \ref{fig:fig3a} that the $(1,1)$ orbit has a period-doubling bifurcation at $\omega \approx 0.066$.)

\vspace{0.1in}

\noindent In contrast as $\omega$ is {\em decreased} then the $(1,3)$ orbit abruptly vanishes and is replaced by the $(1,2)$ periodic orbit at $\omega \approx 0.114$  which is not close to the $(1,3)$ orbit. This is a very dramatic change in the behaviour of the system, and we shall see presently that it is due a
non-smooth {\em grazing} bifurcation as the parameters vary. Again, from computations it appears that the {\em left hand boundary} of the region of existence of the $(1,n)$ orbits is in general due to a grazing bifurcation. We will discuss this point later.

\vspace{0.1in}

\noindent {\em Varying $d$}

\vspace{0.1in}

\noindent Whilst it is analytically convenient to vary the forcing frequency $\omega$ (and also the amplitude $\mu$ of the external forcing), there is no physical evidence of change in the frequency of the external Milankovich forcing over the period of the MPT. Instead we can conjecture that variations in other internal parameters of the PP04 model might have led to this transition. A possible such parameter is $d$ which comes into the definition of the function $F$ which describes the salty bottom water formation efficiency. Accordingly we now fix $\omega = 0.1476$ and $\mu = 0.467$ to be the physically relevant values described in \cite{morupisi2020analysis}. In \cite{paillard2004antarctic} the value of $d = 0.27$ is taken as representative of the current values and at this value we have a stable $(1,3)$ periodic orbit. 
In Figure \ref{fig:fig3a} we plot the resulting Mont\'e-Carlo bifurcation diagram, varying $d$ between 0.2 and 0.4. In this we see a very similar picture to that seen when $\omega$ varies. In particular as series of $(1,n)$ periodic orbits which are suddenly 'created' at grazing bifurcations, and lose stability at period-doubling bifurcations, as $d$ varies. 

\vspace{0.1in}

\noindent A series of numerical experiments shows very similar transitions as other parameters vary. This is due to the {\em generic} nature of the grazing bifurcation which we will explore for the remainder of this paper.

\begin{figure}[htb!]
\centering
\includegraphics[scale=0.4]{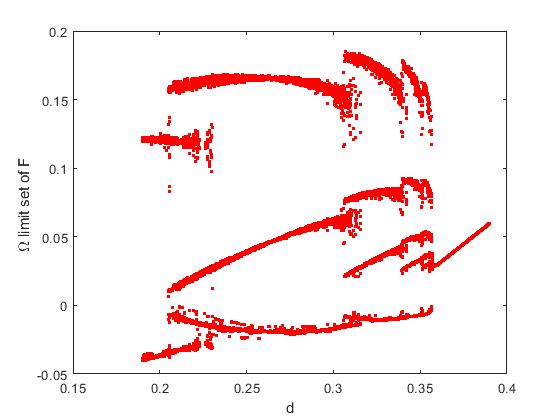}
\caption{Mont\'e-Carlo bifurcation diagram when $\mu=0.467$ and $\omega = 0.1476$ take 'physical' values, and we vary the parameter $d$. In this we see a $(1,2)$ orbit which vanishes at a period-doubling bifurcation at $d \approx 0.23$, a  $(1,3)$ orbit which arises at a {\em grazing} bifurcation at $d \approx 0.2$ and at a period doubling bifurcation at $d \approx 0.31$, followed by windows of existence of $(1,4), (1,5), (1,6)$ orbits. }
\label{fig:fig3a}
\end{figure}

\subsection{Informal explanation of the behaviour at a grazing bifurcation} 

\begin{center}
\begin{figure}[htb!]
\centering
\includegraphics[scale=0.15]{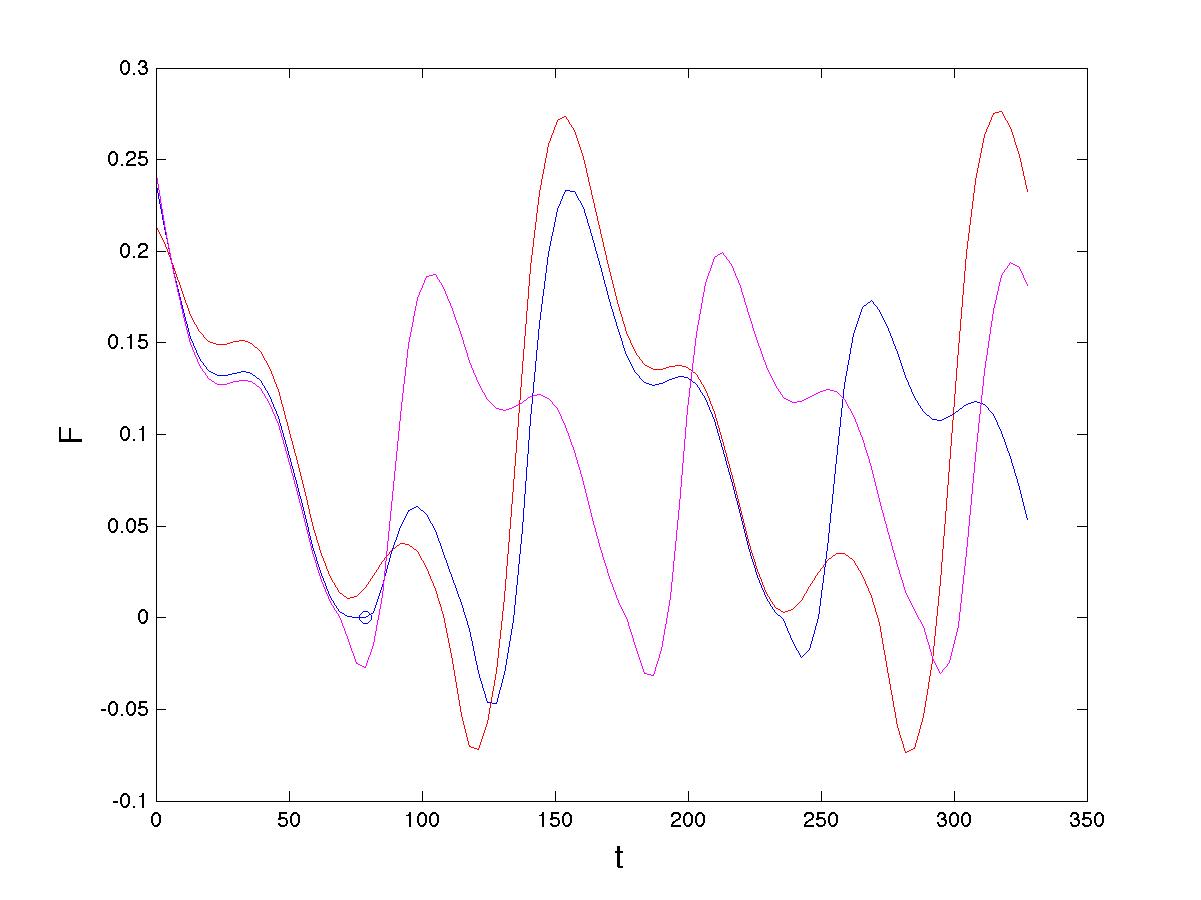}
\caption{Trajectories close to the grazing $(1,3)$ orbit in which we see the large change in behaviour at the grazing value of $V_g = 0.3786.$  }
\label{fig:fig3b}
\end{figure}
\end{center}

\noindent An examination of the form of the $(1,3)$ orbit in Figure \ref{fig:ped2} (right) gives an insight into the cause of the abrupt change in behaviour observed at the grazing bifurcation. In the figure we see that the $(1,3)$ orbit has an impact with the surface $F = 0$ when $t = 200$. For $200 kyr< t < 310 kyr$ it lies in the glacial region with $F > 0$ and for $310kyr < t < 360 kyr$ in the inter-glacial region with $F < 0$. In the glacial region it has a number of oscillations, all of which have $F > 0$, however one of these at $t = t_m = 290kyr$ has a minimum close to $F = 0$. If any parameter (for example $\omega$ or $d$), or the initial state, is varied, then it is likely that there will be a value of the parameter/initial condition and a time $t_g$ close to $t_m$ at which there is a tangential intersection, a {\em grazing event}, when
$F(t_g) = \dot{F}(t_g) = 0.$
If a parameter, or the initial condition, is varied further then the trajectory enters the inter glacial region for $t$ close to $t_g$. It follows from (\ref{Fdisc}) that there is a discontinuity in $\ddot{F}$ at $t_g$ which leads to a significant change in the behaviour of the orbit. To make this more precise we consider the $(1,3)$ orbit in more detail. For $\mu = 0.3$ and $\omega = 0.115$ a direct calculation shows that at $t = 0$ on this orbit we have
$$(V(0),A(0),C(0)) = {\mathbf X}_{3,1}(t=0) = [0.3636,0.2089,0.2356].$$
Now we keep all parameters fixed apart from $V(0)$ and plot the trajectories for $V(0) = 0.3636$ (red) $ V(0) = 0.3786$ (grazing: blue) and $V(0) = 0.38$ (maroon). The resulting trajectories are shown in Figure \ref{fig:fig3b}.  
\noindent In this figure the trajectory at $V \equiv V_g = 0.3785$ has a grazing intersection at $t_g = 78$. The three trajectories, despite starting at very similar values of $V(0)$ have quite different asymptotic behaviour. That starting from $V(0)$ is asymptotic to the $(1,3)$ orbit, the other two are asymptotic to the period $(1,2)$ orbit, but with a different phase in each case.

\subsection{Consequences of the grazing transition}

\noindent Before we develop the theory of the grazing transition in more detail, we now make a series of numerical observations of the dynamical consequences of the grazing events. Examples of this type of behaviour have already been observed, without analysis, in the papers \cite{crucifix2012oscillators}, \cite{paillard2004antarctic}, \cite{mitsui2014dynamics}.

\vspace{0.1in}

\noindent {\em Asymptotic behaviour:} The plot in Figure \ref{fig:fig3b} demonstrates the dramatic changes in behaviour associated with grazing. If we take initial data close to grazing trajectory then the final asymptotic state starting from this initial data can be very different. Thus the grazing orbits mark the separation between the domains of attraction of the different asymptotic states, for example the $(1,2)$ and $(1,3)$ orbits. 

\vspace{0.1in}

\begin{center}
\begin{figure}[htb!]
\centering
\includegraphics[scale=0.12]{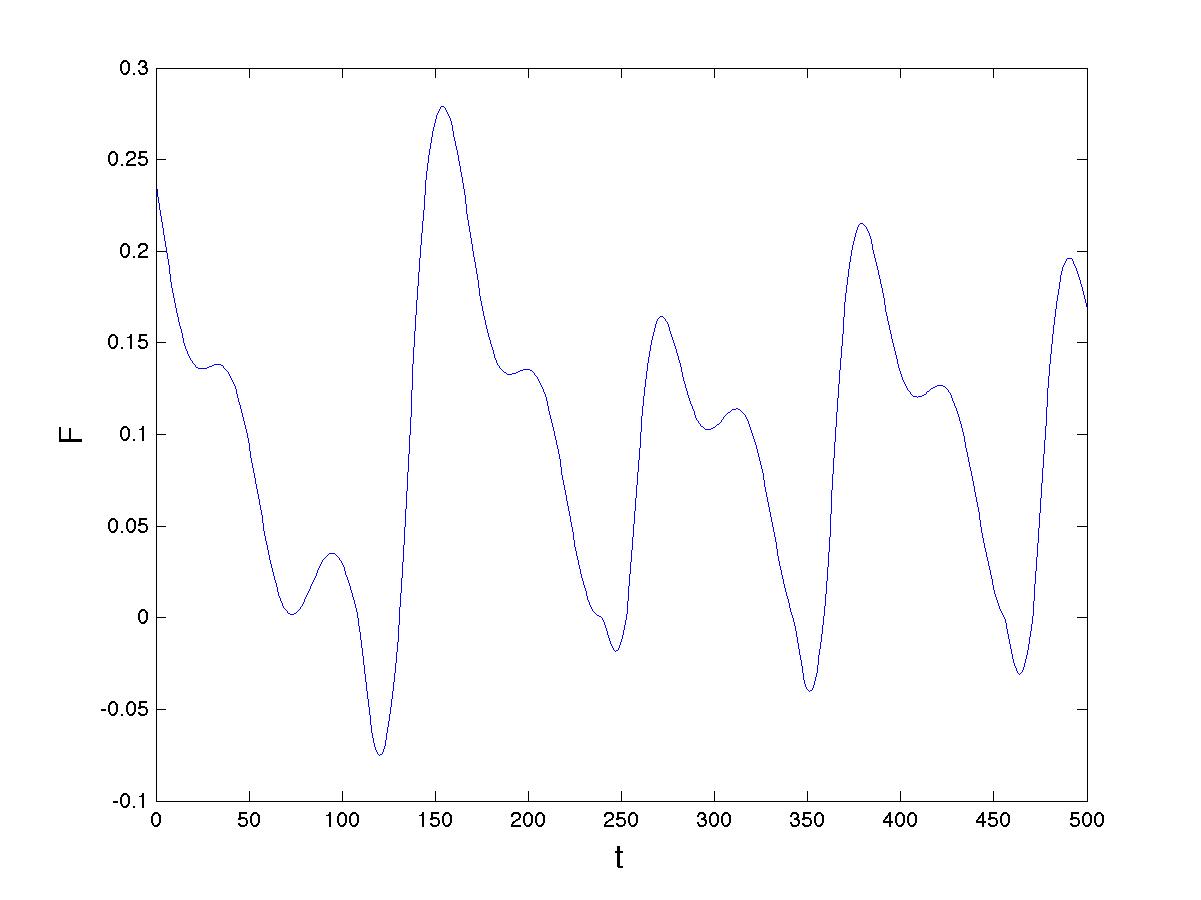}
\includegraphics[scale=0.12]{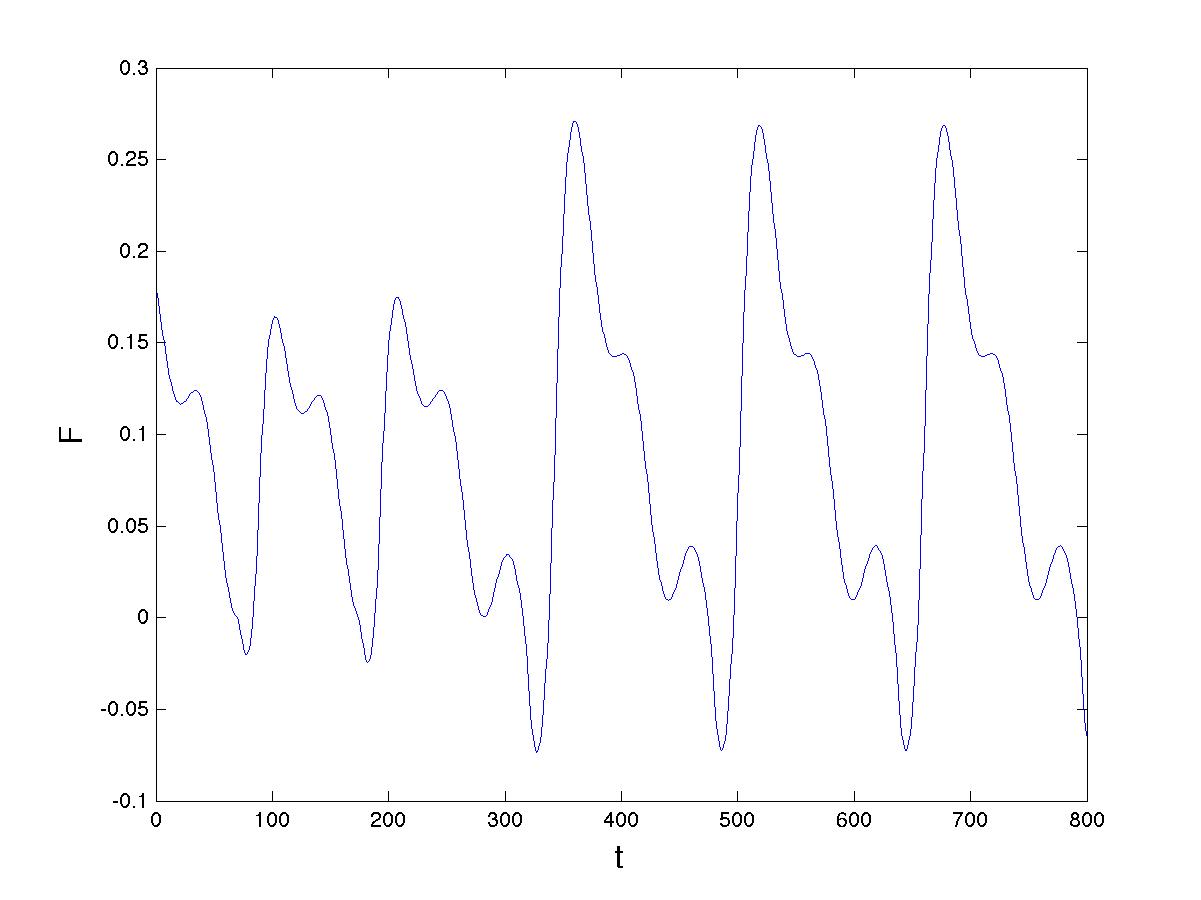}
\caption{Transient changes initiated by grazing events when $\mu = 0.3, d = 0.27$. (Left) $\omega = 0.1138$, (Right) $\omega = 0.119$. }
\label{fig4b}
\end{figure}
\end{center}

\noindent {\em Transient behaviour:} Grazing events can also lead to changes in the transient behaviour of the PP04 system.
We illustrate this in Figure \ref{fig4b} by firstly considering $\omega = 0.1138$, $\mu = 0.3$, $d = 0.27$. At these parameter values the $(1,2)$ orbit is the only stable periodic orbit. If we start with initial data close to the $(1,3)$ orbit for $\omega = 0.115$ we see an initial transient close to the $(1,3)$ orbit, followed by a grazing event, and then a transition to the stable $(1,2)$ orbit.  Similarly we take $\omega = 0.119$ for which only the $(1,3)$ orbit is stable and start close to the $(1,2)$ orbit above. Again we see a transition to the $(1,3)$ orbit triggered by a grazing event. 
 
\vspace{0.1in}

\noindent {\em Slowly varying parameters:} It is reasonable in any climate model to assume that its  parameters can change slowly with time. For example in the Saltzman and Maasch model of the ice ages \cite{saltzman1991first}, there is a slow variable to account for CO2 changes due to tectonic changes, which leads to transitions between states resembling the MPT. As a comparison we introduce slowly varying parameters into the PP04 model. We consider both a change in the 'external' parameters by varying $\omega$ and also a change in the 'internal' parameters by varying $d$.

\vspace{0.1in} 

\begin{center}
\begin{figure}[htb!]
\centering
\includegraphics[scale=0.12]{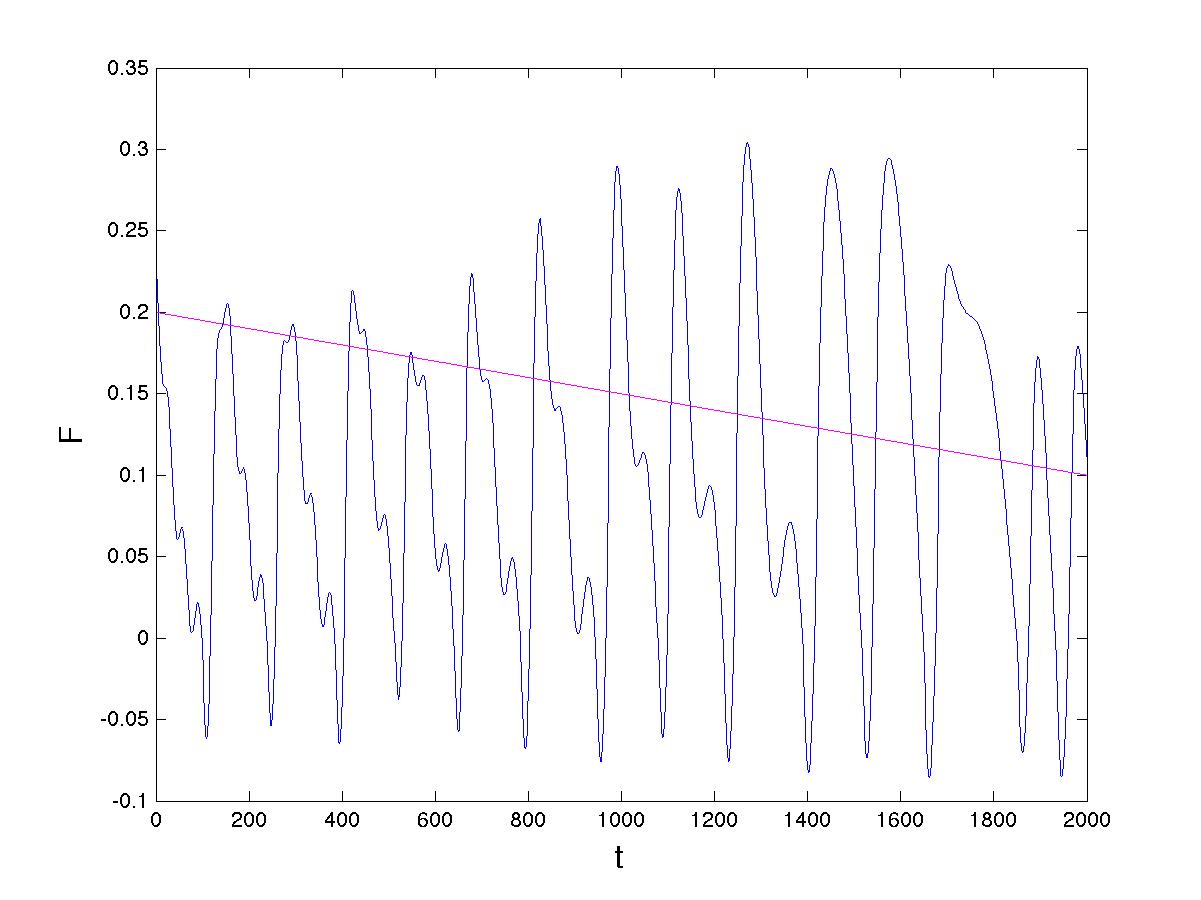}
\includegraphics[scale=0.12]{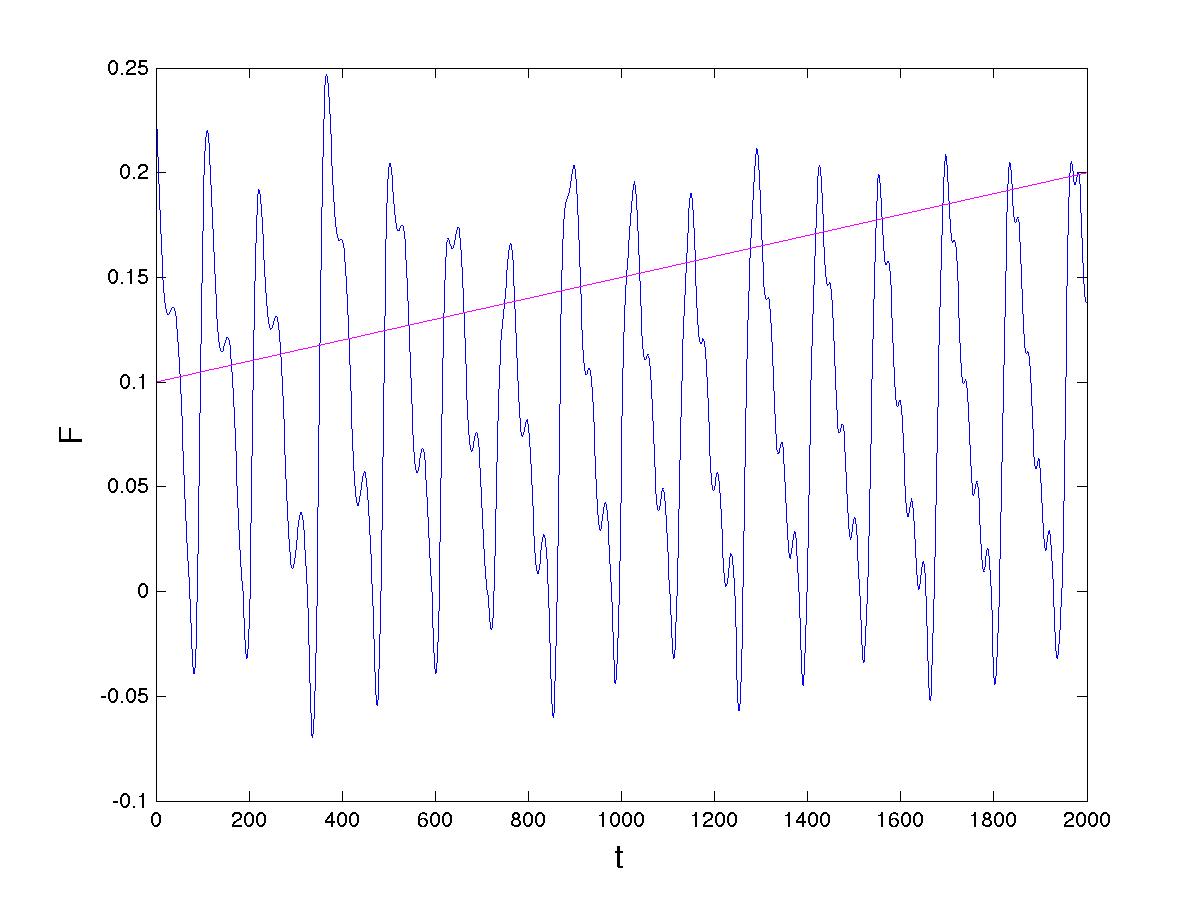}
\caption{The values of $F$ (blue) resulting from solving the PP04 system with a slowly varying $\omega$ (maroon): (left) decreasing, (right) increasing.}
\label{fig4a}
\end{figure}
\end{center}

\noindent {\em Varying $\omega$.} In a first calculation we reduce $\omega$ slowly, setting $\omega = 0.2 - 5*10^{-5} t$ with $0 < t < 2000$kyr. In Figure \ref{fig4a} (left) we show the resulting behaviour of $F$ together with the corresponding value of $\omega$ (maroon).
\noindent In this figure we see an initially stable $(1,3)$ orbit, which transitions
first to a $(1,2)$ orbit at $t \approx 900$kyr ($\omega \approx 0.15$), and then to a $(1,1)$ orbit at $t \approx 1300$kyr ($\omega \approx 0.135$). 
In each case we can see from the figure that these transitions occur just after grazing events. It is interesting to observe that the values of $\omega$ at which the changes occur are higher than the values for the problem with static $\omega$. This is a consequence of the additional dynamics introduced by the variation in $\omega$.

\vspace{0.1in}

\noindent In a second calculation we increase $\omega$ slowly so that
$\omega = 0.1 + 5*10^{-5} t$
again with $0 < t < 2000$kyr. In Figure \ref{fig4a} (right) we show the resulting behaviour of $F$ together with the corresponding value of $\omega$.
 In this figure we see an initially stable $(1,2)$ orbit, which transitions
first to a $(1,3)$ orbit at $t \approx 300$kyr ($\omega \approx 0.115$), and then to a $(1,4)$ orbit at $t \approx 1100$kyr ($\omega \approx 0.155$). 
In each case we can see from the figure that these transitions occur close to grazing events. We observe that the values of $\omega$ at which the changes occur are in this case lower than the values for the problem with static $\omega$. 

\vspace{0.1in}

\begin{figure}[htb!]
\centering
\includegraphics[scale=0.5]{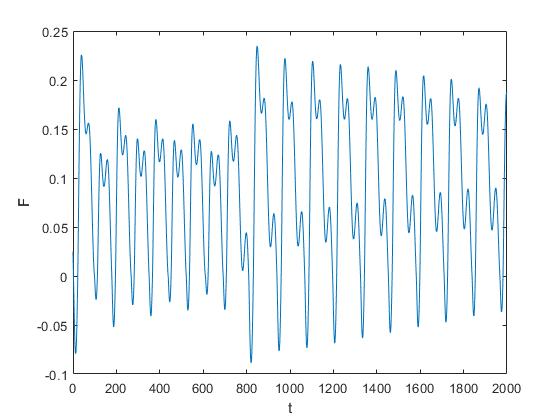}
\caption{The time series solution of $F$ obtained for $\omega=0.1476$ and $\mu=0.467$ with $d$ increasing slowly from $0.2$ to $0.3$. This shows a transition of the solution from a period $(1,2)$ orbit to a period $(1,3)$ orbit at around $t = 800$ thousand years, when $d = 0.24$. The figure has a similar form to the transitions observed at the MPT.}
\label{drange}
\end{figure}

\noindent {\em Varying $d:$} 
Finally we vary an internal parameter by fixing $\omega = 0.1476$, $\mu = 0.467$ and increasing $d$ slowly from $0.2$ to $0.3$ over the range $0 < t < 2000$ (so that $d = 0.2 + 0.1 t/2000$). The results are presented in Figure \ref{drange}. In this figure we see a transition from a low amplitude period $(1,2)$ orbit (or to be more precise a $(2,4)$ orbit which is very close to a $(1,2)$ orbit) to a large amplitude period $(1,3)$ orbit at $t \approx 800, d=0.24$, where there is a near grazing event.  Qualitatively this transition seems to be similar to that observed at the MPT. It is fully consistent with the bifurcation diagram presented in Figure 5, in which we see a $(1,3)$ orbit created at a grazing bifurcation at $d \approx 0.2$, which becomes the global attractor for $0.24 < d < 0.3$. Similar transitions are observed for changes in other system parameters, Whilst we have no direct evidence that changes in $d$ (or any other parameter) are directly responsible for the MPT, the qualitative similarity of the change, and its generic nature, suggest that a grazing bifurcation may have played a role in the change of behaviour observed at the MPT, 

\vspace{0.1in}

\begin{center}
\begin{figure}[htb!]
\centering
\includegraphics[scale=0.15]{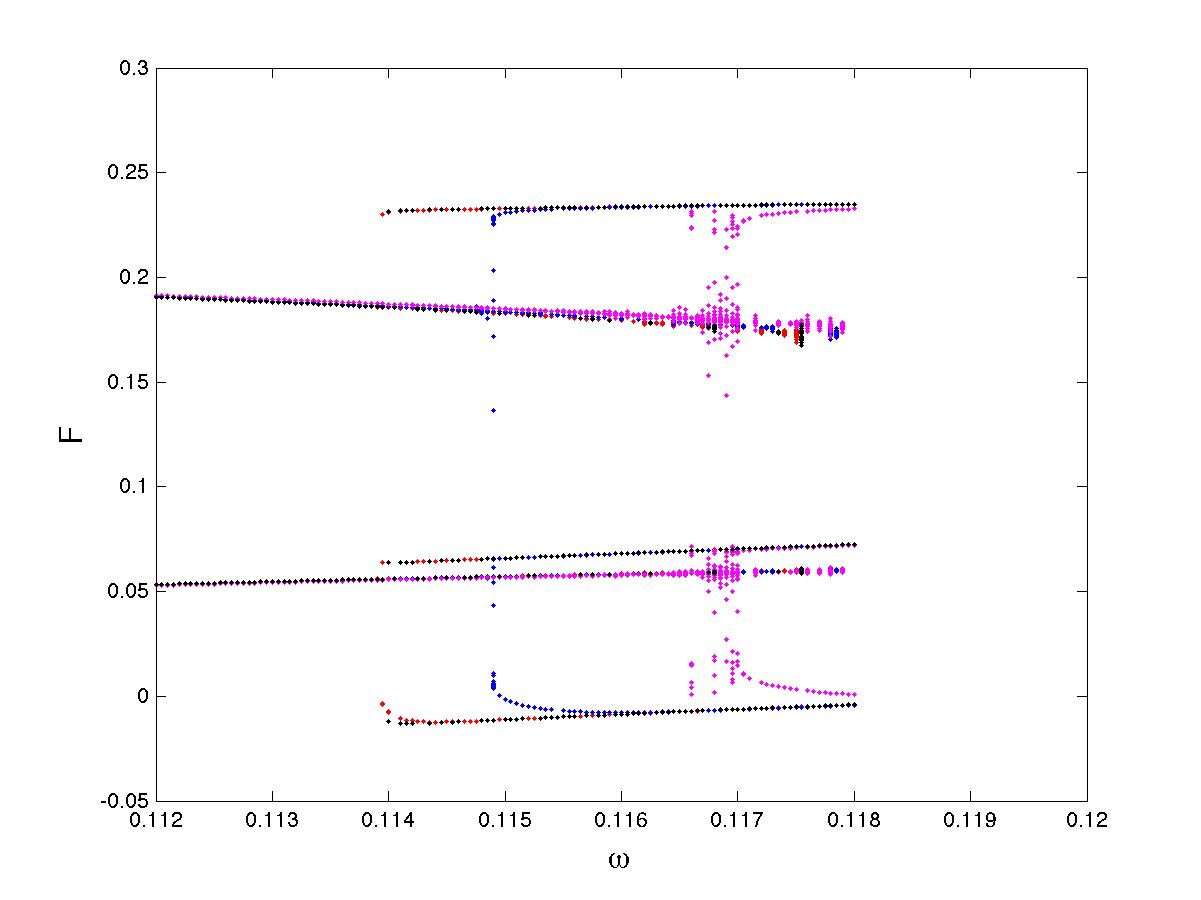}
\caption{Unfolding of the grazing bifurcation. In this figure we see the change as $\eta$ varies from $200$ (maroon), to 500 (blue), 1500 (red), 2000 (black). }
\label{fig:fig4}
\end{figure}
\end{center}

\noindent {\em Tipping points:} There is a significant literature \cite{lenton2008tipping} on climate change associated with {\em tipping points} which are examples of saddle-node bifurcations where for example the eigenvalues of the linearisation of the system around a fixed point pass through zero. Generally before the tipping point some 'warning' is given of a change in the behaviour, due to a 'slowing down' of the system as the eigenvalue approaches zero. In contrast the grazing bifurcation occurs without any warning and is associated with a {\em global} change in the orbit (a grazing intersection) rather than a local change associated with a linearisation. More correctly if we study the smoothed system (\ref{c3aa}) then the grazing bifurcation can be thought of as an unfolding of tipping point, but one in which the behaviour usually associated with tipping occurs in a very small neighbourhood of the bifurcation, which rapidly decreases as $\eta$ increases.   To see a little more of the comparison between the grazing bifurcation and the tipping point, we take a range of values of $\eta$ increasing from $\eta = 200$ to $\eta = 2000$ and plot the resulting Mont\'e-Carlo bifurcation diagram in each case. The result is shown in Figure \ref{fig:fig4}. In this figure we can see that the simple grazing discontinuity when $\eta = 2000$ at $\omega = 0.114$ changes to a much more complex (and harder to analyse) tipping point structure when $\eta = 200$ at $\omega = 0.117$.

\section{The form of the Poincar\'e map close to a grazing event}

\noindent We now make a more formal study of the grazing transition in the PP04 model. The first thing we look at is the precise nature of the change in the system at a grazing event and we discuss this in this section. In the next section we will look at the geometry of the set of initial data which leads to grazing.

\vspace{0.1in}

\noindent As the PP04 model is an example of a discontinuous Filippov system we can make use of existing results on the dynamics of non-smooth systems, given for example in \cite{bernardo2008piecewise}, \cite{filippov2013differential},\cite{simpson2010bifurcations} (although the latter text mainly deals with the case of continuous systems).  Key to our study will be the {\em Poincar\'e Map} $P_S$, which we define as follows. Let a trajectory of the PP04 system be
 ${\mathbf X}(t)$ evolving over a time interval $t \in [t_{\alpha},\;t_{\beta}]$ from an initial point ${\mathbf X}(t_{\alpha}) \equiv {\mathbf X}_{\alpha}$ to a final point 
 ${\mathbf X}(t_{\beta}) \equiv {\mathbf X}_{\beta}$ with 
 $t_{\beta} = t_{\alpha} + 2\pi/\omega$. Such a trajectory will typically intersect $\Sigma$ at a number of intermediate points
$t_{\alpha} < t_i < t_{\omega}$. 

\vspace{0.1in}

\noindent {\bf Definition 4.1} {\em We define the Poincar\'e Map $P_S$ by
\begin{equation}
    P_S {\mathbf X}_{\alpha} = {\mathbf X}_{\beta}
\label{poinc}
\end{equation}
}
\vspace{0.1in}

\noindent We now study the smoothness of $P_S$. Indeed we show that it is a piece-wise smooth map with jump discontinuities associated with  'square-root' behaviour.

\vspace{0.1in}

\noindent {\bf Definition} We define any intersection to be {\em transversal} if $\dot{F}({\mathbf X}(t_i)) \ne 0.$

\vspace{0.1in}

\noindent The following result is well known, see \cite{bernardo2008piecewise}

\vspace{0.2in}
     
\noindent {\bf Lemma 4.2} {\em If there are no intersections of the trajectory ${\mathbf X}(t)$ with $\Sigma$ for $t_{\alpha} < t < t_{\beta}$, or if all intersections are transversal, then the map $P_S:{\mathbf X}_{\alpha} \to {\mathbf X}_{\beta}$ is smooth in a neighbourhood of ${\mathbf X}_{\alpha}$ }

\vspace{0.1in}

\noindent NOTE If there is no impact then the Poincar\'e map is a linear affine map, and from (\ref{H3}) it has the explicit form:
$$P_S {\mathbf X}_{\alpha} = e^{2\pi L/\omega}({\mathbf X}_{\alpha}+L^{-1}{\mathbf b}^{\pm} -{\mathbf c}(t_{\alpha})) 
 -L^{-1} {\mathbf b}^{\pm} + \mathbf{c}(t_{\beta}),$$
 where ${\mathbf c}(t)$ is a smooth function. 

\vspace{0.1in}

\noindent The smoothness of the map $P_S$ changes dramatically at a non transversal grazing intersection. Indeed  we have the following.

\vspace{0.2in}

\noindent {\bf Lemma 4.3} The map $P_S$ is in general discontinuous at those values of ${\mathbf X}_{\alpha} \equiv X_{\alpha,g}$ for which there is a grazing event at time $t_{\alpha} < t_g < t_{\beta}$ at which ${\mathbf X}(t_g) \equiv {\mathbf X}_g.$

\vspace{0.1in}

\noindent {\em Proof} We consider such a trajectory ${\mathbf X}^+$, and associated $F^+$, which we assume lies in $S^+$, with $\ddot{F}^+ > 0.$ Then, from (\ref{Fdisc}), if $\ddot{F}^+(t_g) < r$ it follows that there is an infinitesimally close orbit $X^-$ with $\ddot{F}^-(t)< 0$ which enters the region $S^-$. Suppose that there are no further impacts in the interval $t_g < t < t_{\beta}$ we have
$${\mathbf X}^+(t_{\beta}) = e^{L\Delta}({\mathbf X}_g + L^{-1} {\mathbf b}^+ - {\mathbf c}(t_g)) - L^{-1} {\mathbf b}^+ {\mathbf c} (t_{\beta}),$$
and
$${\mathbf X}^-(t_{\beta}) = e^{L\Delta}(X_g + L^{-1} {\mathbf b}^- - {\mathbf c} (t_g)) - L^{-1} {\mathbf b}^- + {\mathbf c}(t_{\beta}),$$ 
where $\Delta = t_{\beta} - t_g.$ Hence
\begin{equation}
{\mathbf X}^+(t_{\beta})-{\mathbf X}^-(t_{\beta}) = L^{-1} (e^{L\Delta}-I)({\mathbf b} ^+ - {\mathbf b}^-).
\end{equation}
Provided that $\Delta > 0$, then as ${\mathbf b}^+ \ne {\mathbf b}^-$ the map is discontinuous. 

\vspace{0.1in}

\noindent In the more general case, if there is a further impact for either ${\mathbf X}^-$ or ${\mathbf X}^+$, then by the above reasoning, the two orbits will still drift apart before the impact, and there is no reason to expect that they will come together later.   \qed 

\vspace{0.1in}

\noindent NOTE as $L$ has only negative eigenvalues $0 > \lambda_1 > \lambda_2 > \lambda_3$, it follows that if $\lambda_1 \; \Delta  \gg 1$ then
\begin{equation}
{\mathbf X}^+(t_{\beta})-{\mathbf X}^-(t_{\beta}) \approx L^{-1} ({\mathbf b}^- - {\mathbf b}^+).
\end{equation}

\vspace{0.1in}

\begin{center}
\begin{figure}[htb!]
\centering
\includegraphics[scale=0.15]{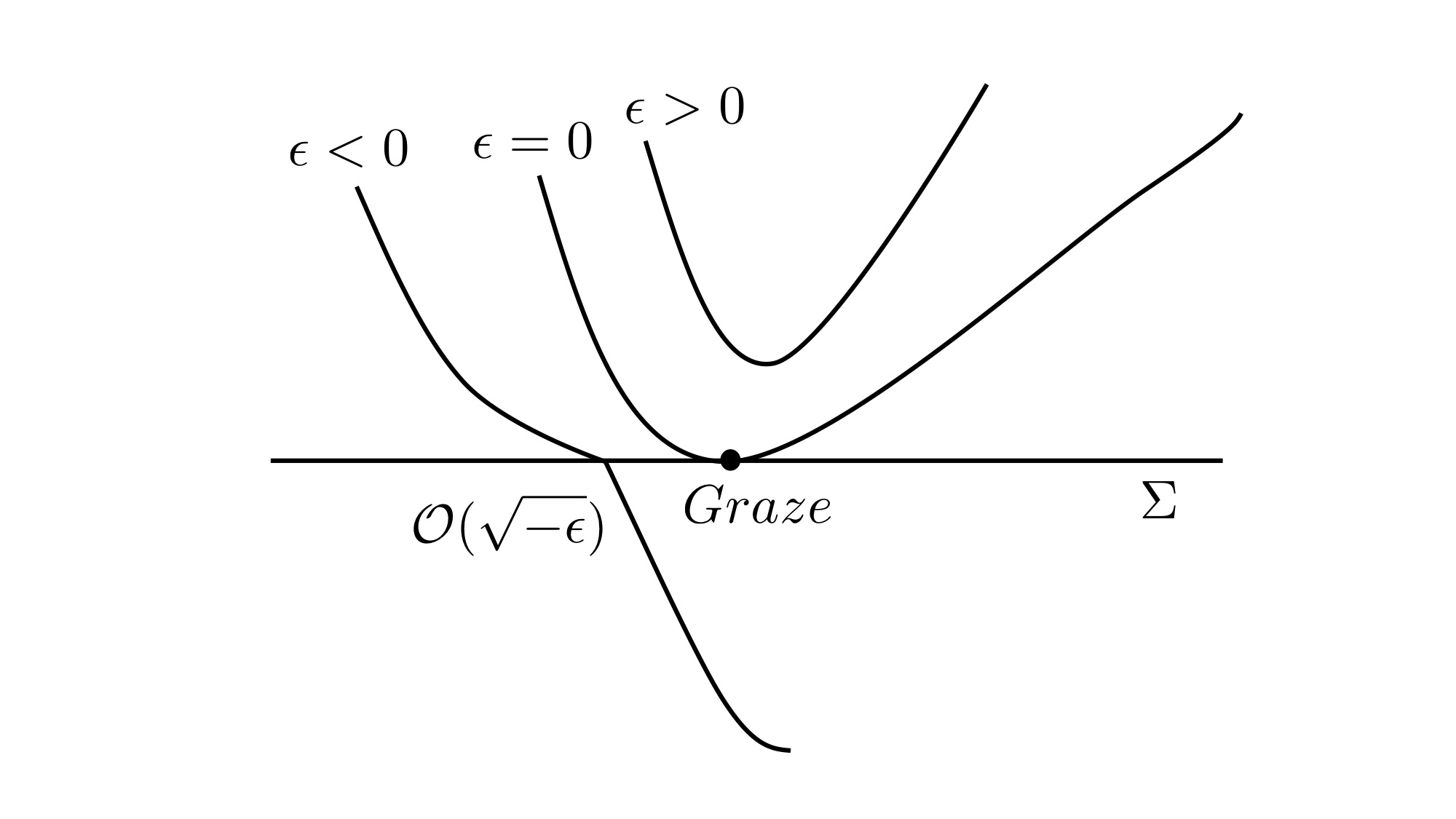}
\caption{A sketch of three trajectories close to the grazing orbit. In this we can see that if $\epsilon>0$ then there is no impact and if $\epsilon < 0$ then there is an impact at a point  ${\cal O}(\sqrt{-\epsilon})$ distant from the grazing point.}
\label{fig:fig3e}
\end{figure}
\end{center}

\noindent We now consider the form of the Poincar\'e map close to such a point ${\mathbf X}_{\alpha}$ in a little more detail. To do so we consider the case where ${\mathbf p}$ is a (one-dimensional) vector and we take as initial data 
$${\mathbf X}_{\alpha} = {\mathbf X}_{\alpha,g} + \epsilon {\mathbf p}.$$
Without loss of generality we will assume that if $\epsilon > 0$ then there is no impact, and if $\epsilon < 0$ then there is an impact at a time $t^*$  slightly before $t_g$. As there is a quadratic minimum at $t_g$ it follows from the construction that there is a $C$ so that to leading order $t^* = t_g - C \sqrt{-\epsilon}$. 
This is illustrated in Figure \ref{fig:fig3e}. We can see from this figure that when $\epsilon > 0$ then the (blue) orbit is close to the (red) grazing orbit. As there is no impact, the map $P_s$ for $\epsilon > 0$ is smooth and hence has a linearisation about the grazing orbit. In contrast if $\epsilon < 0$ then the map has, as described above, a discontinuity. Furthermore the map will include terms of the form $t_{\beta} - t^* = t_{\beta} - t_g + C\sqrt{-\epsilon}.$ Hence there are contributions to the map $P_S$ involving $\sqrt{-\epsilon}$. 
\noindent We conclude from this informal argument that there are vectors ${\mathbf a}$, ${\mathbf b}$, ${\mathbf c}$ and ${\mathbf d}$  so that to leading order
\begin{equation}
    P_S ({\mathbf X}_{\alpha,g}+ \epsilon \; {\mathbf p}) = {\mathbf a} + \epsilon \;  \; {\mathbf b} \quad \mbox{if} \quad \epsilon > 0,
\end{equation}
and
\begin{equation}
 P_S ({\mathbf X}_{\alpha,g}+ \epsilon \;  {\mathbf p}) = {\mathbf c}  + \sqrt{-\epsilon} \; \; {\mathbf d} \quad \mbox{if} \quad \epsilon < 0,    
\end{equation}
where ${\mathbf a} \ne {\mathbf c}.$ A detailed proof of this result is given in \cite{bernardo2008piecewise}, \cite{nordmark1997universal}.

\vspace{0.1in}

\begin{center}
\begin{figure}[htb!]
\centering
\includegraphics[scale=0.15]{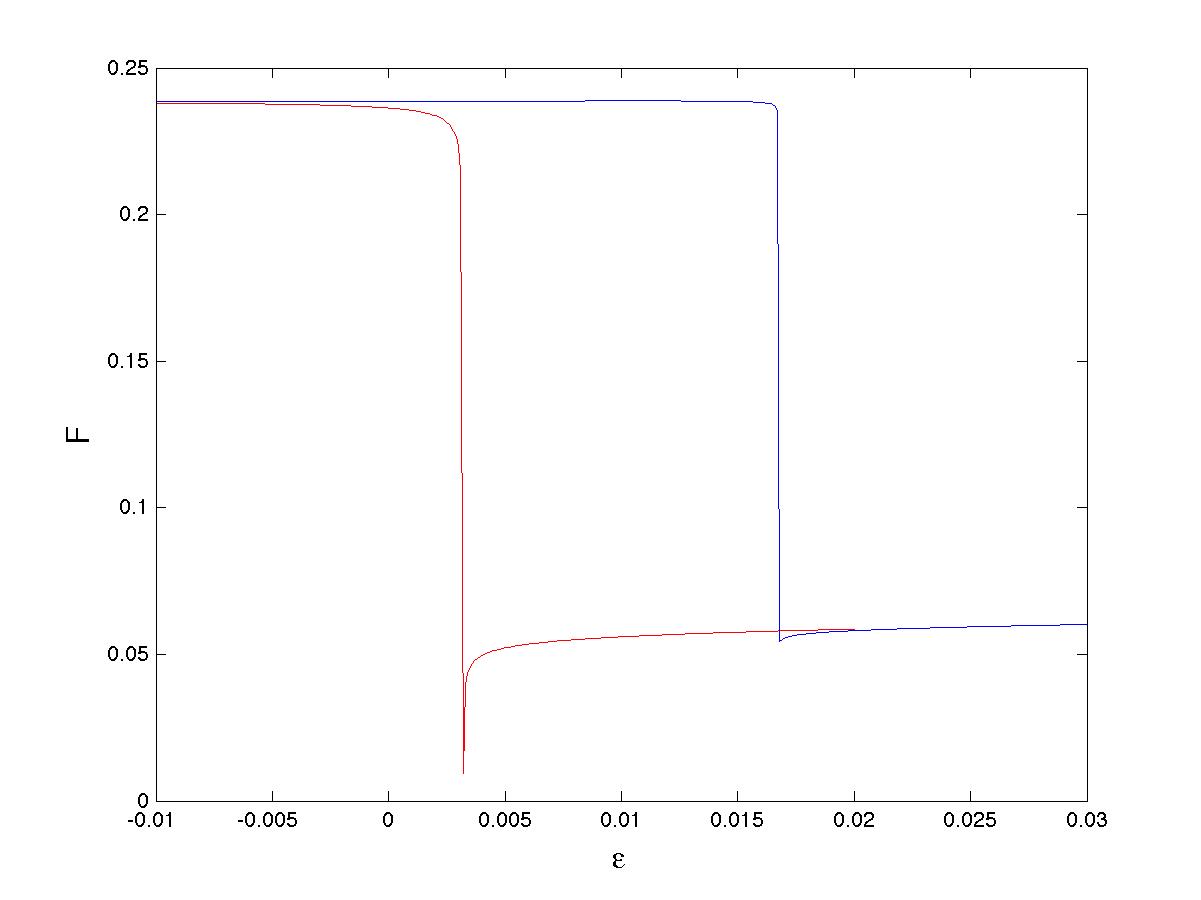}
\caption{The Poincar\'e map close to the grazing event when $\eta = 1500$ (blue) and when $\eta = 500$ (red).}
\label{fig:fig3f}
\end{figure}
\end{center}

\noindent As an example calculation we take the $(1,3)$ orbit above with $t_{\alpha} = 0$ and let ${\mathbf p} = (1,0,0)^T$. The values of $F$ at $t_{\beta}$ are then plotted in Figure \ref{fig:fig3f} in the two cases of $\eta = 1500$ (blue) and $\eta = 500$ (red)

\section{Grazing surfaces}

\noindent We have seen from the above that large changes in the behaviour of the solutions to the PP04 model (both in transient and asymptotic dynamics) close to grazing trajectories which have a grazing event at the time $t_g$. If we fix an initial time $t_0$ then the initial data which leads to such grazing trajectories forms the {\em Grazing Surface} ${\cal G}$.  This surface then plays a critical role in understanding the dynamics of the PP04 model. For example, if, as we vary a parameter, then (the initial data for) a periodic orbit intersects ${\cal G}$, then that orbit will become unstable, and the dynamics of the solution will (as we have seen) change dramatically. It is thus important to determine the geometry of ${\cal G}$. In this section we will establish the following, rather unexpected, result.

\vspace{0.2in}

\noindent {\bf Theorem 5.1} {\em The set ${\cal G}$ is a multi-leaved surface. Each of the leaves is locally planar. On each leaf the grazing time $t_g$ is close to being constant}.

\vspace{0.2in}

\noindent NOTE The feature of the PP04 model which leads to this form of the grazing surface, is that the matrix $L$ is dissipative, and it has widely separated eigenvalues.

\vspace{0.1in}

\noindent We prove Theorem 5.1 through the following sub-sections. 

\subsection{General theory}
\noindent If $t_0$ is fixed, consider a Poincar\'e section given by 
$$ \Pi = \{ (V,A,C):  t = t_0 \}. $$
We now construct the grazing surface ${\cal G} \subset \Pi$ defined as follows
$${\cal G} = \{ (V,A,C) \in \Pi: \exists \; t_g > t_0:  F(t_g) = \dot{F}(t_g) = 0.$$

\noindent NOTE This definition allows for multiple intersections of ${\mathbf X}$ with $\Sigma$ before the grazing impact at time $t_g$.

\vspace{0.2in}

\noindent {\bf Lemma 5.2} {\em The surface ${\cal G}$ is locally well-defined and smooth locally planar 2D sub-manifold of $\Pi$ apart from isolated points at which $\ddot{F}(t_g) = 0$ or if there is a non-transverse impact before $t_g$. }

\vspace{0.2in}

\noindent {\bf Proof} We define 
$${\mathbf G}(V,A,C,t_g) \equiv (F(t_g), \dot{F}(t_g)).$$
We now perturb ${\mathbf X} = (V,A,C)^T$ by  $\delta {\mathbf X} \equiv (\delta V, \delta A, \delta C)$ and assume that an impact occurs at a time $t_g + \delta t$.
Suppose that $t_g$ is the time of the first non-transverse impact. Then $F_{\mathbf X}$ is a well defined operator. As $dF/dt(t_g) =0$ it follows that
$$\delta {\mathbf  G} = (F_{\mathbf X} \delta {\mathbf X}  \quad \dot{F}_{\mathbf X} \delta {\mathbf X} + \ddot{F} \delta t) + {\cal O}(\|\delta {\mathbf X} \|^2 + \delta t^2).$$
If $\delta {\mathbf  G} = {\mathbf 0}$ then $\delta {\mathbf X}$ must lie locally on the plane normal to the vector $F_{\mathbf X}$ so that $F_{\mathbf X} \delta {\mathbf X} = 0.$ The value of $\delta t$ is then given by
$$\delta t = \dot{F}_{\mathbf X}/\ddot{F}.$$
The surface ${\cal G}$ is thus well defined, provided that $\ddot{F} \ne 0.$

\vspace{0.1in}

\noindent The construction fails if $\ddot{F} = 0.$ Similarly the construction fails if there is a non-transverse impact before $t_g$. As we have seen already, in this case $F_{\mathbf X}$ is no longer well defined and the analysis above breaks down. \qed

\vspace{0.2in}

\noindent {\bf Lemma 5.3} {\em Let ${\cal G}_0$ be the subset of ${\cal G}$ for which there are no impacts in the interval $t_0 < t < t_g$. Then ${\cal G}_0$ is close to being a collection of planes, and the grazing impact time $t_g$ is close to a constant on each plane.}

\vspace{0.2in}

\noindent {\em Proof} The general solution of the PP04 model  between impacts is given by
\begin{equation}
    {\mathbf X}(t) = e^{L(t-t_0)}\left({\mathbf X}_0 - {\mathbf f}(t_0) \right) + {\mathbf f}(t) 
\end{equation}
where ${\mathbf f}(t) = -L^{-1} {\mathbf b}^{\pm} + {\mathbf c}(t)$ is an explicitly computable function which does not depend on ${\mathbf X}_0$. It follows that
\begin{equation}
    F(t) = c^T e^{L(t-t_0)}({\mathbf X}_0 - {\mathbf f}(t_0) ) + c^T{\mathbf f}(t) + d,
\end{equation}
and
\begin{equation}
    \dot{F}(t) = c^T L e^{L(t-t_0)}({\mathbf X}_0 - {\mathbf f}(t_0) )  + c^T \dot{\mathbf f}(t)
\end{equation}
Now we can set
$$L = U^{-1} \lambda U \quad \mbox{where} \quad \Lambda = diag(-\lambda_1,-\lambda_2,-\lambda_3).$$
\noindent It follows from the choice of constants used in the PP04 model that we have
$$0 < \lambda_1 \ll \lambda_2 \ll \lambda_3.$$
As a consequence, for even moderately large values of $t$ we have that
$$e^{L(t-t_0)} \approx e^{-\lambda_1(t-t_0)} U \; diag(1,0,0) \; U^{-1} \equiv e^{-\lambda (t - t_0)} M.$$
\noindent Similarly
$$e^{L(t-t_0)} L \approx -\lambda_1 e^{-\lambda (t - t_0)} A$$
where 
$$M = U \; diag(1,0,0) \; U^{-1}.$$
Hence, if
$${\mathbf n} = {\mathbf c}^T M$$
then the condition for a graze at the time $t_g$ is given (approximately)
by the expressions
$$
    e^{-\lambda_1(t_g - t_0)}  {\mathbf n}^T ({\mathbf X}_0 - {\mathbf f}(t_0)) + {\mathbf c}^T {\mathbf f}(t_g) + d= 0,$$
    and  
    \begin{equation}
    -\lambda_1 e^{-\lambda_1(t_g - t_0)}  {\mathbf n}^T ({\mathbf X}_0 - {\mathbf f}(t_0) )  + {\mathbf c}^T \dot {\mathbf f}(t_g) = 0.
    \label{oc1}
\end{equation}

\noindent It is clear that (\ref{oc1}) defines two parallel planes, orthogonal to ${\mathbf n}$ in which
${\mathbf X}_0$ must lie. As these planes have the same normal vector ${\mathbf n}$ to coincide it follows that $t_g$ must satisfy the equation:
\begin{equation}
 -\lambda_1 \left( {\mathbf c}^T {\mathbf f}(t_g) + d \right) =  {\mathbf c}^T \dot {\mathbf f}(t_g).  
 \label{oc2}
\end{equation}

\noindent The equation (\ref{oc2}) gives a leading order approximation for $t_g$. There will of course be small variations in $t_g$ due to the effect of the smaller exponential terms which have been set to zero above.  The equation (\ref{oc2})  will have a number of discrete solutions for $t_g$. Each such solution 'defines' a plane of solutions ${\mathbf X}$ and  $t_g$ is approximately constant on each such plane. \qed

\vspace{0.2in}

\noindent For the default values given for the PP04 model we have the structure 
\begin{equation}
M = [*  | {\mathbf 0}  | *] \quad \mbox{and} \quad {\mathbf n} = (-0.478, 0, 0.228)^T.
\label{nvals}
\end{equation}

\vspace{0.1in}

\noindent This is a significant result as it shows that {\em each 'plane' of the grazing surface is to leading order independent of the values of the antarctic ice $A$.
} This is a key aspect of the geometry of ${\cal G}$. 

\vspace{0.1in}

\noindent  Significantly the result that the grazing time $t_g$ is almost constant on each 'leaf' of ${\cal G}$ applies even if there are several impacts on the trajectories starting from that leaf. 

\vspace{0.2in}

\noindent {\bf Lemma 5.4} If the time of the last impact before grazing is $t_{\alpha}$ and if $t_g - t_{\alpha}$ is not small, then to leading order $t_g$ is a solution of the equation (\ref{oc2}).

\vspace{0.1in}

\noindent {\em Proof} Let $t_{\alpha}$ be either $t_0$ if there is no impact, or the last impact time before $t_g$. Similarly let ${\mathbf X}_{\alpha}$ be ${\mathbf X}(t_{\alpha})$ and let $\Delta = t_g - t_{\alpha}$ we have
$${\mathbf X} = e^{L \Delta} ({\mathbf X}_{\alpha} - {\mathbf f}(t_{\alpha})) + {\mathbf f}(t).$$
Thus, provided that $\Delta$ is sufficiently large so that 
$\exp(-\lambda_2 \Delta) \ll \exp(-\lambda_1 \Delta)$
it follows that
$$F(t) = c^T e^{L\Delta} ({\mathbf X}_{\alpha} - {\mathbf f} (t_{\alpha})) + c^T {\mathbf f} + d \approx e^{-\lambda_1} n^T ({\mathbf X}_{\alpha} - {\mathbf f}(t_{\alpha})) + c^T {\mathbf f}(t) + d $$
and similarly
$$\dot{F}(t) \approx -\lambda_1 e^{-\lambda_1 \Delta} n^T( X_{\alpha} - {\mathbf f}(t_{\alpha})) + c^T \dot{{\mathbf f}}.$$
As $F(t_g) = \dot{F}(t_g) = 0$ the result follows.

\qed

\vspace{0.1in}

\noindent We note that the nonlinear equation (\ref{oc2}) will have multiple solutions, which leads to the 'leaf' structure of ${\cal G}$. Note further that as the function ${\mathbf f}(t)$ is periodic in $t$ then if $t_g$ is a solution, then so is $t_g + 2\pi k/\omega$ for integer $k$. This observation leads immediately to the following definition.

\vspace{0.1in}

\begin{figure}[htb!]
\centering
\includegraphics[scale=0.15]{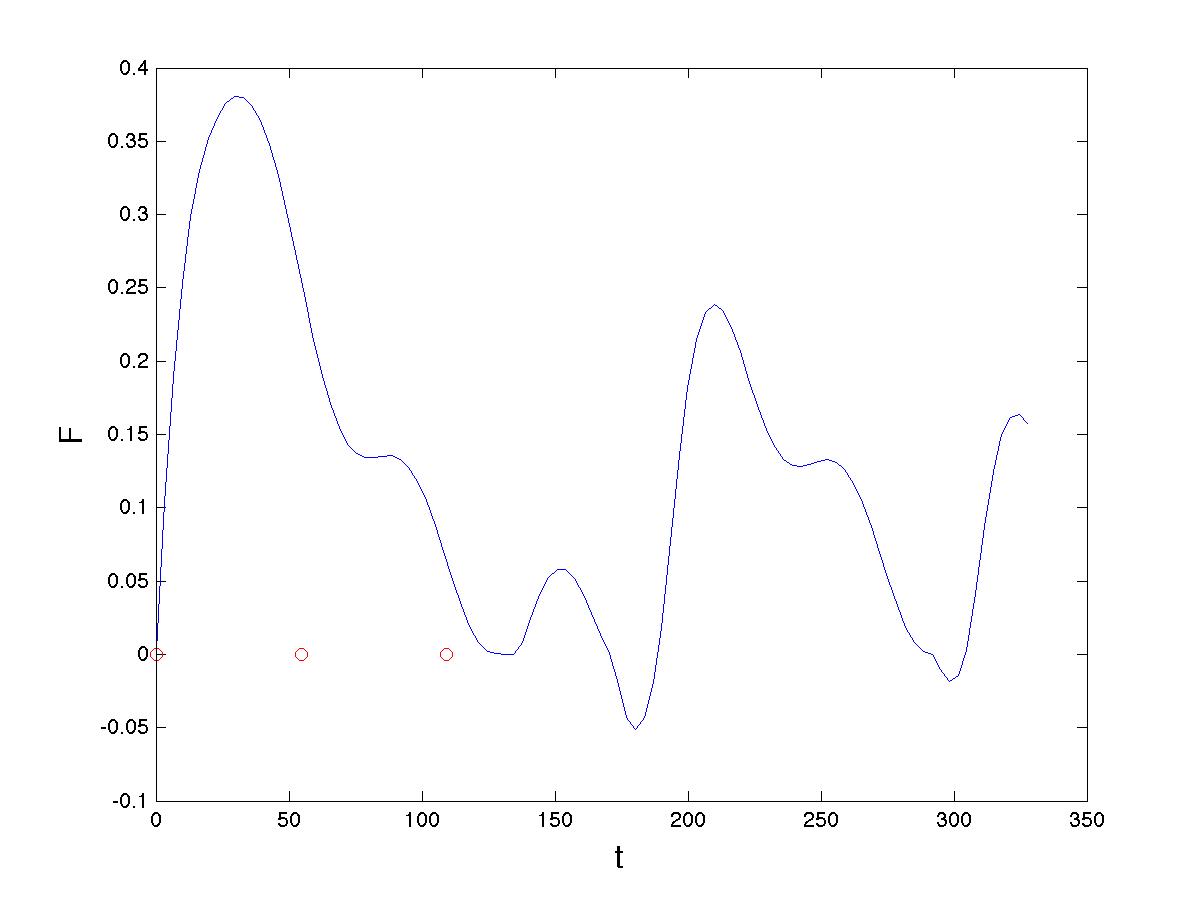}
\caption{A grazing orbit starting from $t=0$ with a graze at the time $t=125$. Orbits starting at times $2\pi/\omega$ and at time $4\pi/\omega$ also lead to grazes at the same time.}
\label{fig:gorbit}
\end{figure}

\noindent {\em Push backs and pull forwards:} Consider any trajectory ${\mathbf X}(t)$  starting from  a leaf of the set ${\cal G}$ at the time $t_0$ so that ${\mathbf X}_0 \equiv {\mathbf X}(t_0) \in {\cal G}$. It follows immediately that any point 
${\mathbf X}_{-k} = P^{-k} {\mathbf X}_0 \in {\cal G}$. Of course the new grazing time will then be the time $t_g + 2\pi k/{\omega}$. Similarly
${\mathbf X}_{k} = P^{k} {\mathbf X}_0 \in {\cal G}$ provided that $t_0 + 2\pi k/\omega < t_g.$ We call the resulting leaves of ${\cal G}$ the respective pull-backs and push forwards of the original leaf.   We deduce that ${\cal G}$ comprises a set of disconnected surfaces which are ${\cal G}_0$ together with its various pull backs and push forwards. We note further from the previous discussion that the surfaces comprising ${\cal G}_0$ will be close to being planes, and each will be orthogonal to ${\mathbf n}$. The push forwards will also all be planes as the map $P_S$ in this case is linear. However the pull backs may not be, as there may be additional impacts, which introduce nonlinearity into $P_S$.  Figure \ref{fig:gorbit} together with the times $t = 2\pi k /\omega$ to illustrate the pull back and push forward points.
\vspace{0.1in}

\noindent Finally we compute $t_g$ more accurately.

\vspace{0.1in}

\noindent {\bf Corollary 5.5} {\em Suppose that there are constant vectors ${\mathbf p,q,r}$ such that
$$ {\mathbf f}(t) = {\mathbf p} \cos(\omega t) + {\mathbf q} \sin(\omega t) + {\mathbf r},$$
then to leading order  $t_g$ satisfies the trigonometric identity:
\begin{equation}
    {\mathbf c}^T ( \lambda_1 \; {\mathbf p} + \omega \; {\mathbf q}) \cos(\omega t) + {\mathbf c}^T (\lambda_1  \; {\mathbf q} - \omega \; {\mathbf p}) \sin(\omega t) + \lambda_1 \;  ({\mathbf c}^T {\mathbf r} + d ) = 0.
    \label{cnov1}
\end{equation}
}

\vspace{0.1in}

\noindent {\bf Proof} This follows immediately from substitution of ${\mathbf f}$  in (\ref{oc2}).

\section{Numerical computation of the grazing surface}

\noindent We now calculate the grazing surface ${\cal G}$ numerically for various parameter values. To do this we fix $t_0 = 0$ and compute the surface leaves ${\mathbf X}_{0} =(V(t_0),A(t_0), C(t_0))$, and associated times $t_g$ so that for each point on the  leaf we have $F(t_g) = {\dot{F}}(t_g) = 0.$ 
As we have shown above, ${\cal G}_0$ we expect that $t_g$ will be lose to constant on each leaf.  For representing the surface it is convenient to take a section of each set determined as above. For convenience we take the section given by $A = A_0$, for some fixed $A_0$. As we have shown above, the resulting surface will only depend weakly upon the value of $A$.  In general the intersection of each leaf with this section will be a curve $(V,C)$ which is locally close to being a line, with $t_g$ close to being a constant on the line.

\vspace{0.1in}

\noindent As an example calculation  we take $t_0 = 0$, $\mu = 0.3$, $\omega = 0.115$.
For these values we can solve equation (\ref{cnov1}) directly and we find that the solutions include
$$t_g \in \{17.777, 40.6454, 72.4141, 95.2818, 127.0505, 149.4812, 181.6869, 236.3233 \}.$$

\noindent For the above values of $\omega$ and $\mu$, as we have seen above, there is a $(1,3)$ periodic solution which has a near grazing intersection, and a $(1,2)$ periodic solution which is not close to grazing. The $(1,3)$ solution has three intersections ${\mathbf X}_{3,i} = (V,A,C)^T$, $i=1,2,3$ with the surface $t_0=0$ given by 
 $${\mathbf X}_{3,1} = (0.3636,0.2089,0.2356) \in S^+,$$
 $${\mathbf X}_{3,2} = (0.7042, 0.5934, 0.04747) \in S^+$$ 
 $${\mathbf X}_{3,3} = (0.811, 0.749, 0.1786) \in S^-.$$
Correspondingly, the $(1,2)$ periodic solution has two intersections given by $X_{2,1} = [0.4025,  0.2972,  0.214] \in S^+$ and $X_{2,2} = [0.718,  0.611,  0.0398] \in S^+.$ If we now fix $A = 0.2089$ then, as $X_{3,1}$ leads to an orbit with a near grazing impact, we expect that the intersection of the grazing surface ${\cal G}$ with the surface $A = 0.2089$ will contain a curve which lies close to $X_1$. Keeping $C = 0.2356$ fixed we find that there are  several values of $V \equiv V_i$ leading to grazing intersections at times $t_{g,i}$. The computed values of these are given by 
$$(V_1,t_{g,1}) = (0.3786, 72.4),  \quad  (V_2,t_{g,2}) = (0.97,22), $$
$$(V_3,t_{g,3}) = (-0.7,127.051),  \quad  (V_4,t_{g,4}) = (0.525,181.687),$$
$$(V_5,t_{g,5}) = (0.925,181.687), \quad  (V_6,t_{g,6}) = (-0.39,236.3233).$$
The  orbits from the last four points, all have at least one non transversal impact before $t_{g,i}$.

\vspace{0.1in}

\begin{figure}[htb!]
\centering
\includegraphics[scale=0.2]{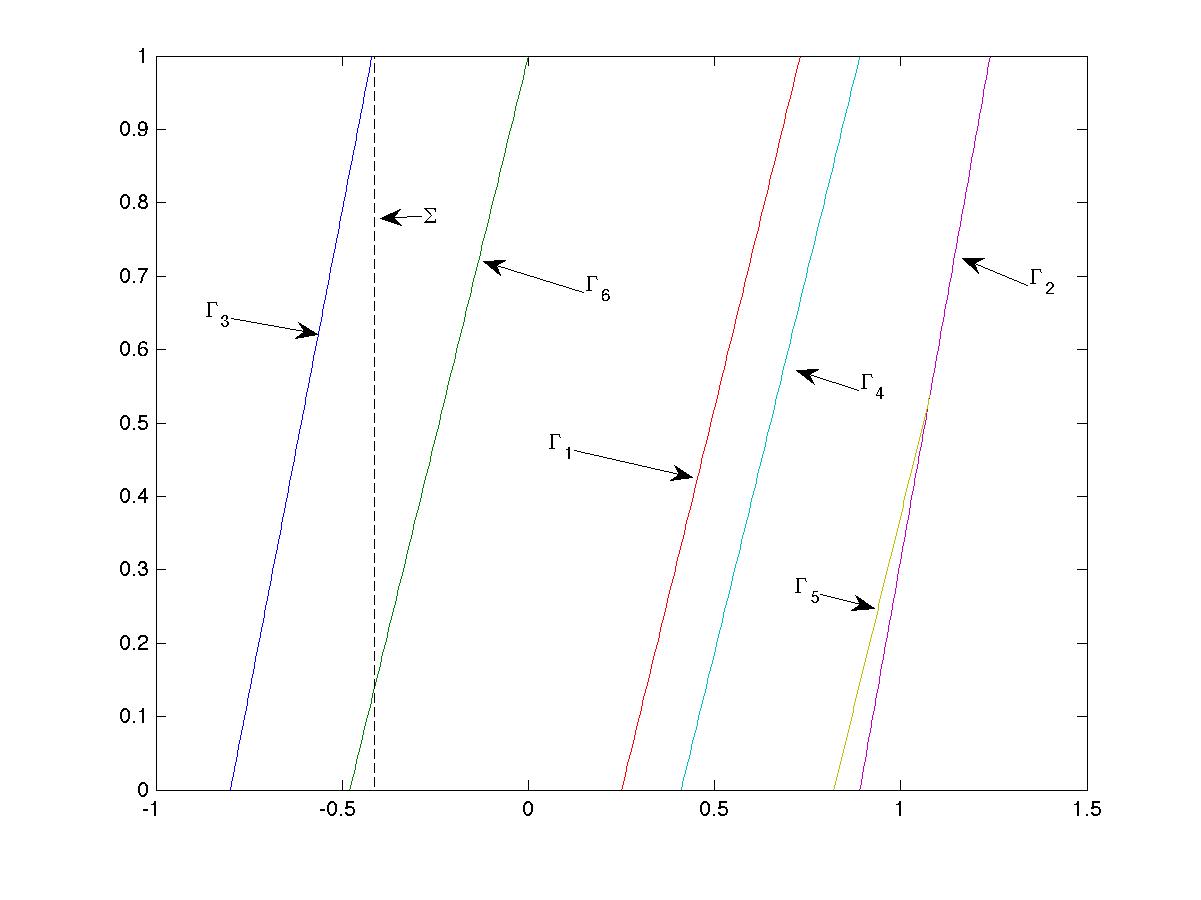}
\caption{A section through the set {\cal G} at $A = 0.2089$ showing the 
the six curves $\Gamma_i$, and the normal vector ${\mathbf n}$ and $\Sigma$. Note the strong linearity of each curve, and the intersection of $\Gamma_2$ and $\Gamma_5$.}
\label{fig:gline}
\end{figure}

\noindent Letting $C$ vary, but keeping $A = 0.2089$ fixed, we now calculate the 6 leaves of ${\cal G}$ containing each of the points $V_i$ above.  In Figure \ref{fig:gline} we show the intersection of ${\cal G}$ with the surface $A = 0.2089$. As expected, this surface comprises a series of curves $\Gamma_i$ containing the points $(V_i,C_i)$.  We note the strong linearity of each of the curves in this Figure, even for those that involve an impact before $t_{g,i}$. The normal vector ${\mathbf n}$ calculated above is, as expected, orthogonal to each line, apart from $\Gamma_2$. This is to be expected, as this line has the shortest time interval before $t_{g,2}$ and hence will be the worst approximated by the results obtained by assuming that the time interval between impacts is large.  For completeness we also plot the intersection of $\Sigma$ with the surface $A =0.2089$ in this case. 

\vspace{0.1in}

\noindent We observe that the curves $\Gamma_2$ and $\Gamma_5$ intersect.  At this point the curve $\Gamma_5$ ceases to exist. This is consistent with the earlier theory. What has happened at the intersection point of $G_2$ and $G_5$ is that an orbit starting from this set of initial conditions has a non-transversal at $t=t_{g,2}$ before the second non-transversal impact at $t = t_{g,5}$. As we will see this impacts in an interesting way on the shape of the domains of attraction of the periodic orbits.

\vspace{0.1in}

\noindent Numerical experiments with a variety of different values of the parameters lead to very similar pictures for ${\cal G}$ with linear curves, some of which intersect in a similar manner to that seen above when $\Gamma_5 \cap \Gamma_2$.

\section{Domains of attraction (DOA)}

\noindent We now consider in more detail the impact of grazing on the asymptotic behaviour of the solutions of the PP04 model as we vary either the initial data. We have already seen in Section 4 that small changes in the initial data close to the grazing orbit can lead, asymptotically, to either a change in the type of periodic orbit, or a change in the phase of the periodic orbit. To illustrate this behaviour we consider the domains of attraction of the different types of periodic orbit. In particular we take a  fixed time $t_0$ and we consider those values of the initial data $(V_0,A_0,C_0)$ which lead, asymptotically, to trajectories which have $\omega-$limit sets each of the different types of periodic orbit (not differentiating between the same periodic orbit but with different phase.)
We show that the form of the domain of attraction is very closely linked to the structure of the grazing set ${\cal G}$.

\vspace{0.1in}

\begin{figure}[htb!]
\centering
\includegraphics[scale=0.2]{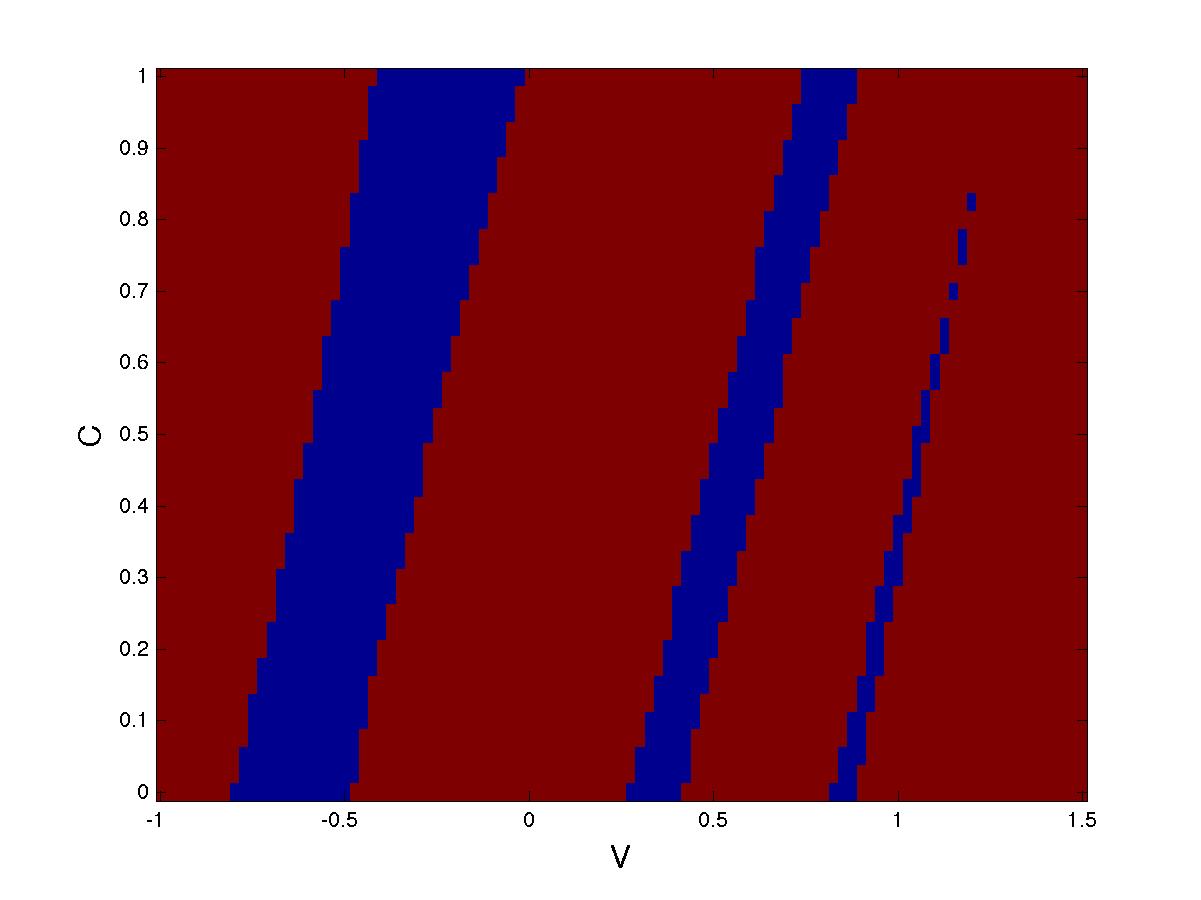}
\caption{The domains of attraction of the $(1,3)$ orbit (red) and the $(1,2)$ orbit (blue). It can be seen that the boundaries of these domains correspond exactly to the curves $\Gamma_i$ shown in Figure \ref{fig:gline}. In this calculation, A = 0.2089 $\mu = 0.3$ $\omega = 0.115$ $\eta = 1500$} 
\label{fig:15001}
\end{figure}

\noindent As a first calculation we consider the earlier example and compute the domains of attraction of the periodic orbits, when $\omega = 0.115$ and $\mu=0.3$. For ease of calculation we fix $A = A_0$ and divide the region of initial data $(V_0,C_0)$ into small squares. We then evolve the trajectory forward in time starting from the centre point of each square, and continuing the calculation until the trajectory has converged to its asymptotic state. We then apply a {\tt fft} to this state, and determine the type of the final orbit by looking at the resulting peaks in the spectrum.  The resulting domains of attraction of the $(1,2)$ (blue) and the $(1,3)$ (red) orbits, restricted to the plane $A = 0.2089$ and taking $t_0 = 0$, are shown in Figure \ref{fig:15001}. We see a sharp linear separation between these domains, the boundaries of which  exactly coincide with the curves $\Gamma_i$ in the set ${\cal G}$, including the 'V-shaped' section of the DOA when $V \approx 1$. Interestingly, there is no evidence in Figure \ref{fig:15001} of the discontinuity set $\Sigma.$ This is because transversal intersections with $\Sigma$ do not have any significant influence on the long term dynamics. 

\vspace{0.1in}

\begin{figure}[htb!]
\centering
\includegraphics[scale=0.1]{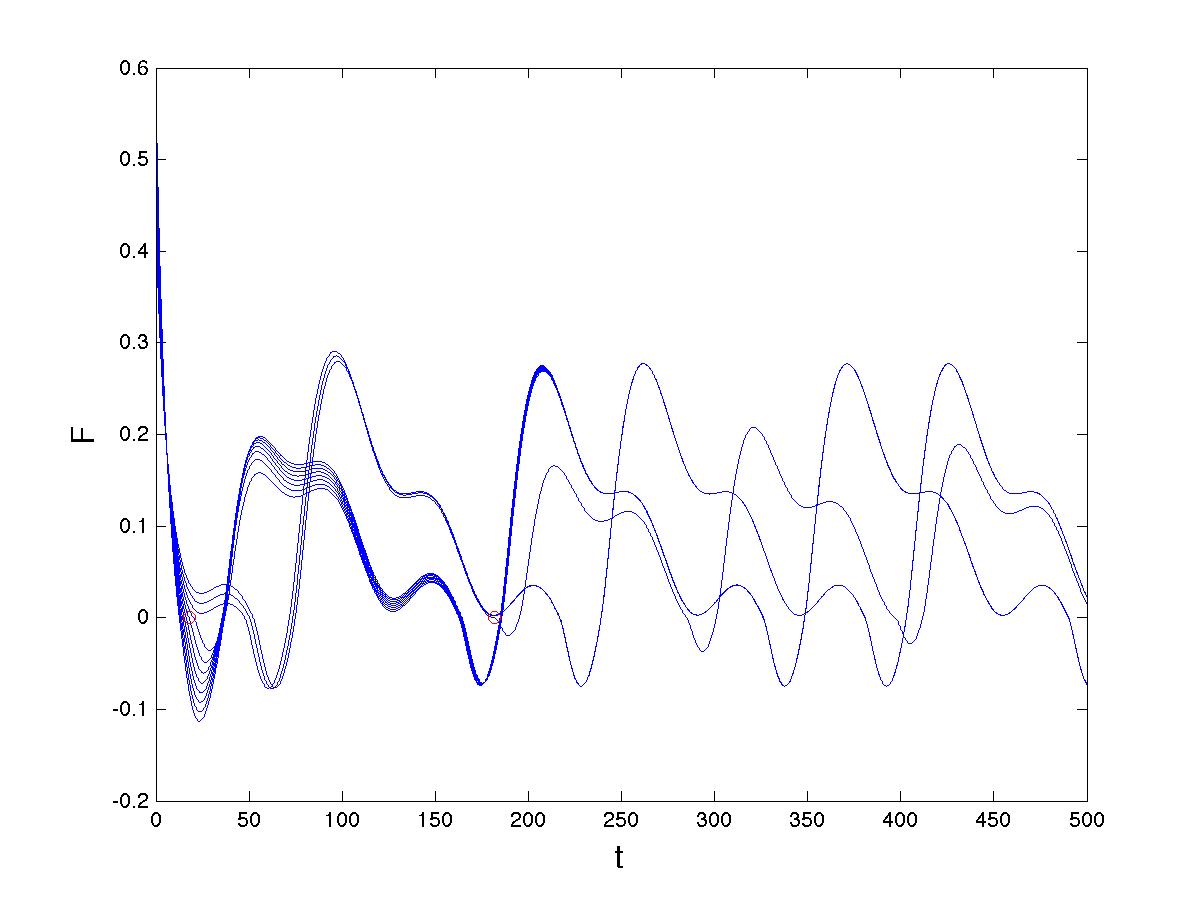}
\includegraphics[scale=0.1]{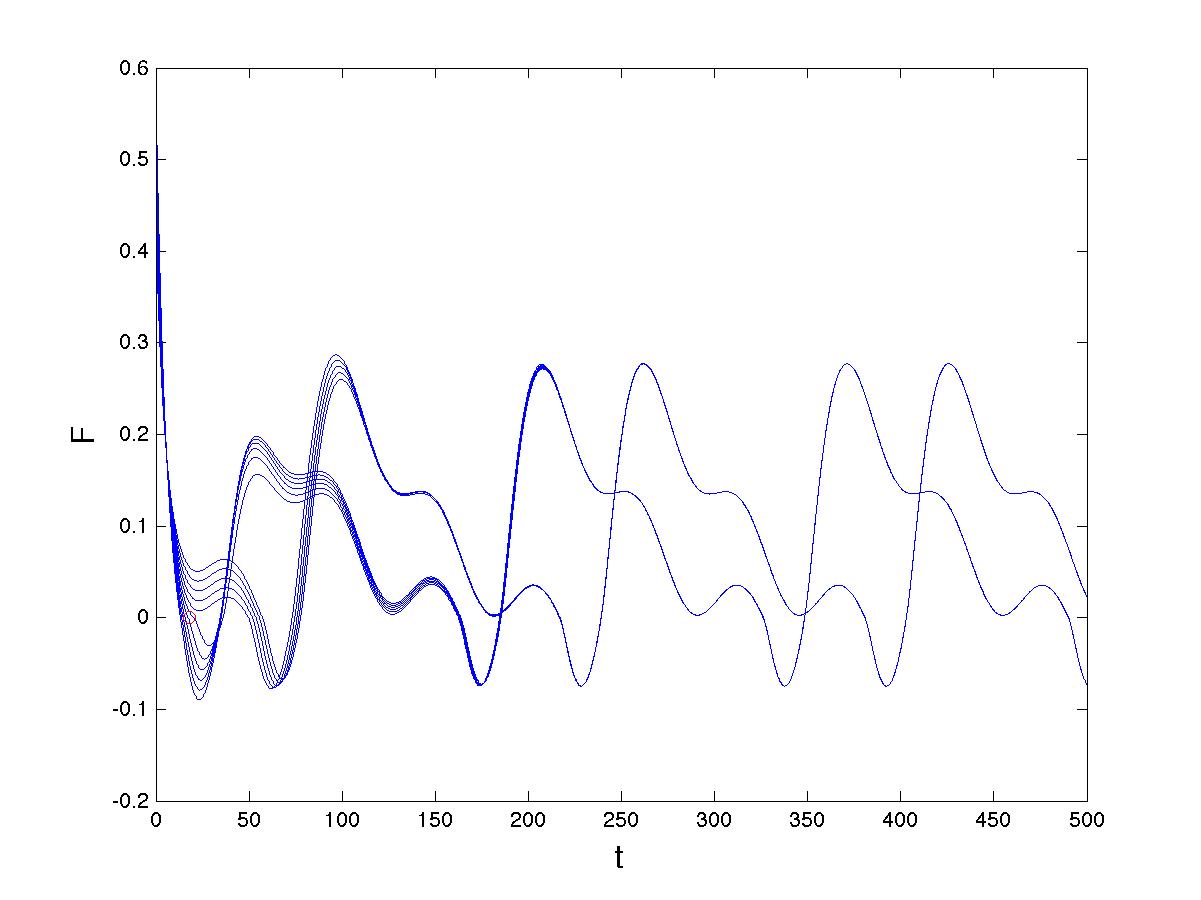}
\caption{A cross section through the Domain of attraction with $C = 0.6$ (left) and $C = 0.9$ (right). In each case we take a range of values of $1 < V < 1.5$. These show transitions between the $(1,3)$ and the $(1,2)$ orbits as well as between orbits of different phase.} 
\label{compo}
\end{figure}

\noindent For most of this figure we see that the lines which are the intersection of the grazing set ${\cal G}$ with this section, mark a change from the $(1,3)$ to the $(1,2)$ orbit or vice-versa. An orbit starting close to the boundary of the DOA (and hence close to ${\cal G}$) will exhibit transient dynamics associated with one orbit before asymptotic dynamics associated with the other. Such changes in the evolutionary behaviour then resembles those that occur at the MPT. This does not seem to be the case close to the rightmost such grazing line (say close to $V = C = 1$). However, what we do see here is a sudden change in the {\em phase} of the $(1,3)$ orbit. To illustrate this in Figure \ref{compo} we show two cross-sections through this figure taking a range of values of $V$ varying between $1 < V < 1.5$, and considering the two cases of $C = 0.6$ and $C= 0.9$. If $C=0.6$ we see a transition from $(1,3)$ to $(1,2)$ to $(1,3)$ with a different phase. If $C=0.9$ we see the change of phase, but no $(1,2)$ orbit.  In both cases we show the calculated grazing times.

\begin{figure}[htb!]
\centering
\includegraphics[width=0.7\textwidth]{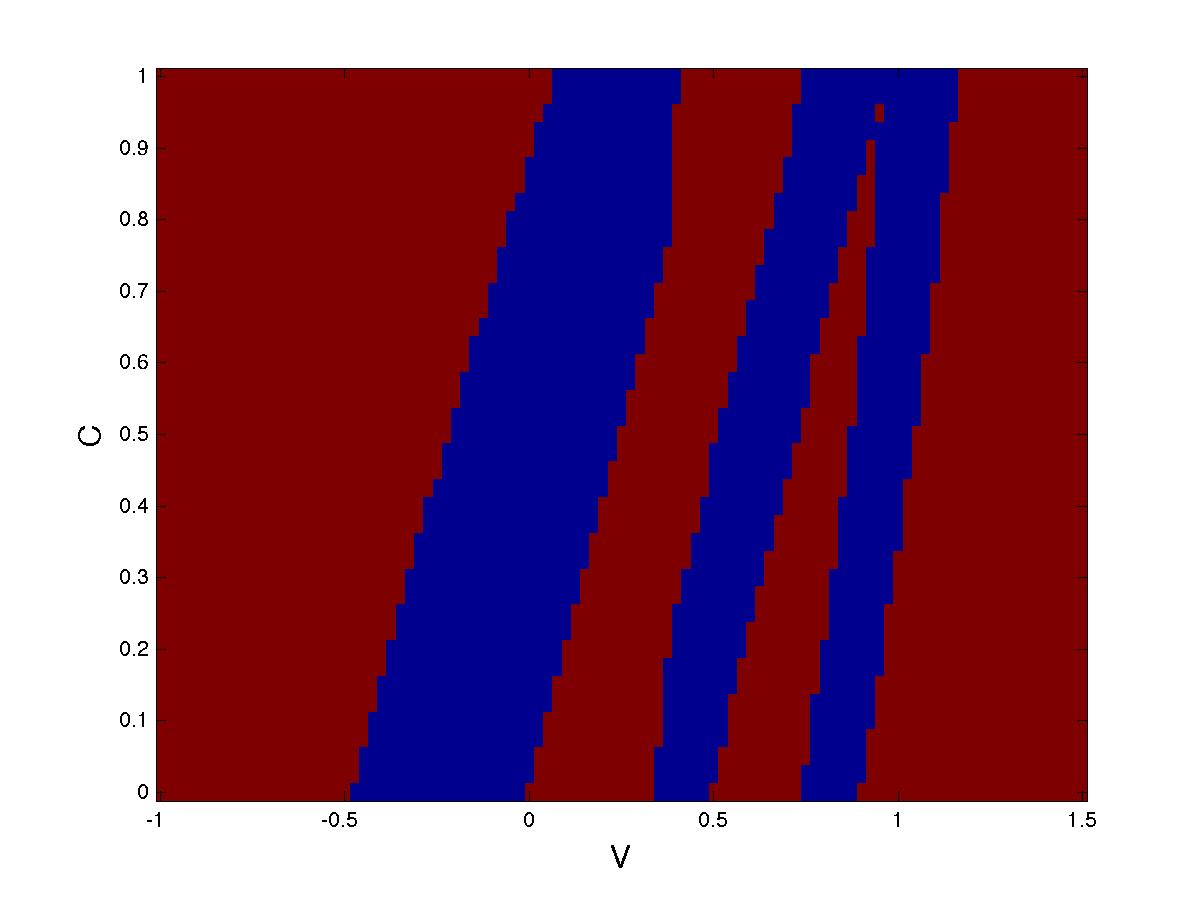}
\caption{The domains of attraction of the $(1,3)$ orbit (red) and the $(1,2)$ orbit (blue). In this calculation we take the 'physically relevant' values of, A = 0.55 $\mu = 0.467$ $\omega = 0.124$ $\eta = 1500$} 
\label{fig:phys1}
\end{figure}

\noindent For comparison we show the DOA for the case of $\mu = 0.467$ and $\omega = 0.124$, which are the physically relevant values considered \cite{morupisi2020analysis} for which $A = 0.55$ is used, see  Figure\ref{fig:phys1}. Again we see that the boundaries between the different domains of attraction are  approximated by straight lines orthogonal to ${\mathbf n}$, which correspond to the surfaces $\Gamma_i$ of ${\cal G}$. As before we see an intersection of these lines, which leads to a V-shaped region in the DOA.

\section{'Arnold Tongues' and curves of grazing bifurcations}

\begin{figure}[htb!]
\centering
\includegraphics[scale=0.28]{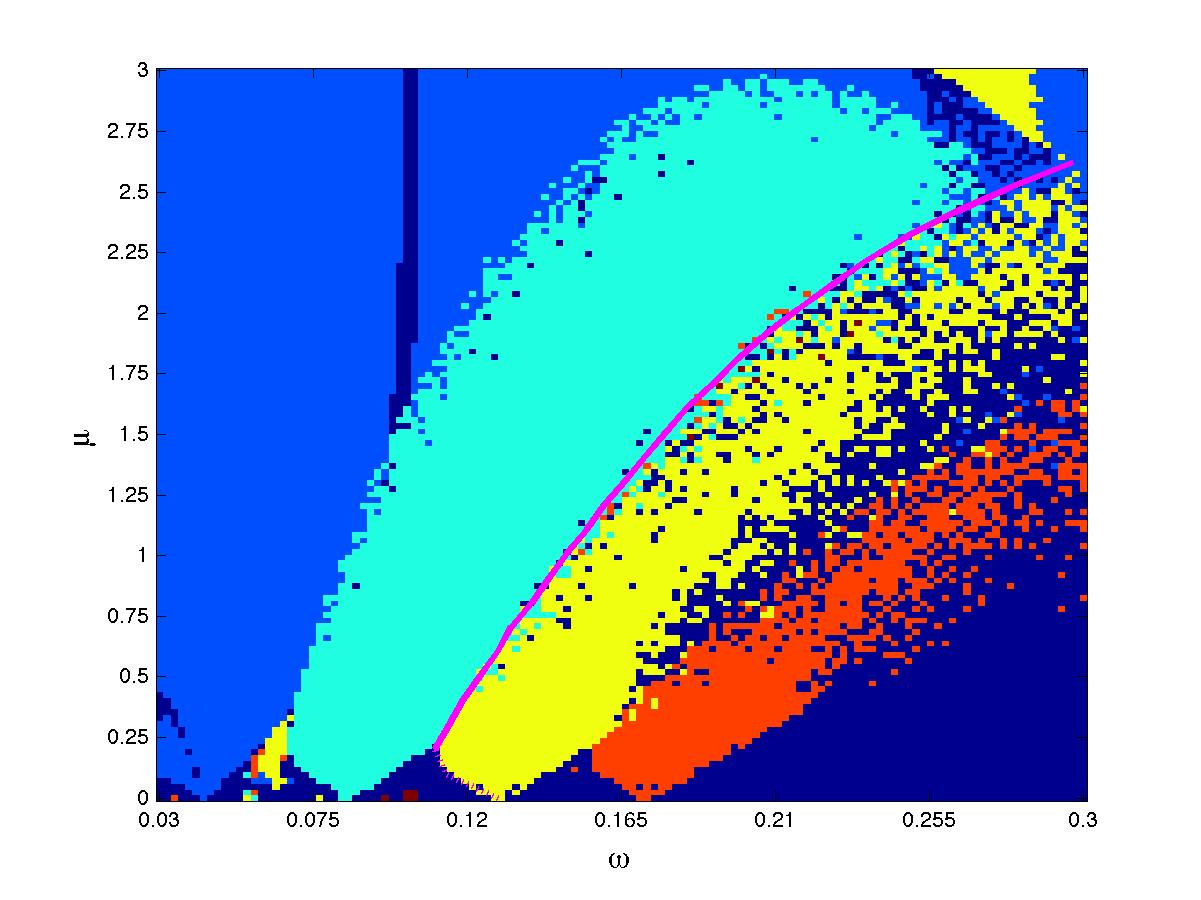}
\caption{Tongues of existence of the $(1,1),(1,2),(1,3),(1,4)$ and $(2,5)$ orbits, showing the 'leaf-like' structure of the tongues for the $(1,2)$, $(1,3)$ and $(1,4)$ orbits, the large region of existence of the $(1,1)$ orbit and the small region of existence for the $(2,5)$ orbit. The grazing curve $G_3$ indicated  by the solid mauve line marks the left-hand boundary of the region of existence of the $(1,3)$ orbit for the $\mu$ approximately greater than 0.25. The dotted mauve line shows the curve of saddle-node bifurcations.} 
\label{fig:leaf}
\end{figure}

\noindent We next consider the change in the behaviour that arises as we vary the parameters in the PP04 model. Again grazing bifurcations play a critical role in this. We again consider varying $\omega$ and $\mu$ as the main parameters, and look again at the regions of existence of the various periodic $(m,n)$ orbits. A more complete picture of these regions is shown in Figure \ref{fig:leaf}. 

\vspace{0.1in}

\noindent For small values of the forcing $\mu$ the tongues of existence of the $(m,n)$ orbits lie close to the 'resonant' values of $n \omega^*/m$ with regions of existence bounded by lines of smooth bifurcations (saddle-node and period doubling). This is the 'standard' picture that arises from an Arnold Tongue calculation. See \cite{schilder2007computing,morupisi2020analysis} for example.

\vspace{0.1in}

\noindent However, as $\mu$ increases we see a change from the 'standard picture'. There is a first value of $\mu \approx 0.1$ at which the leftmost boundaries of each of the existence regions for the solutions become determined by grazing bifurcations on the periodic orbit. We now explore this situation in more detail.

\vspace{0.1in}

\noindent {\bf Definition 8.1} The curve $G_n = \{ (\omega,\mu): \omega = G_n(\mu) \}$ is the set of points for which the $(1,n)$ orbit has an additional grazing impact.

\vspace{0.1in}

\begin{figure}[htb!]
\centering
\includegraphics[scale=0.25]{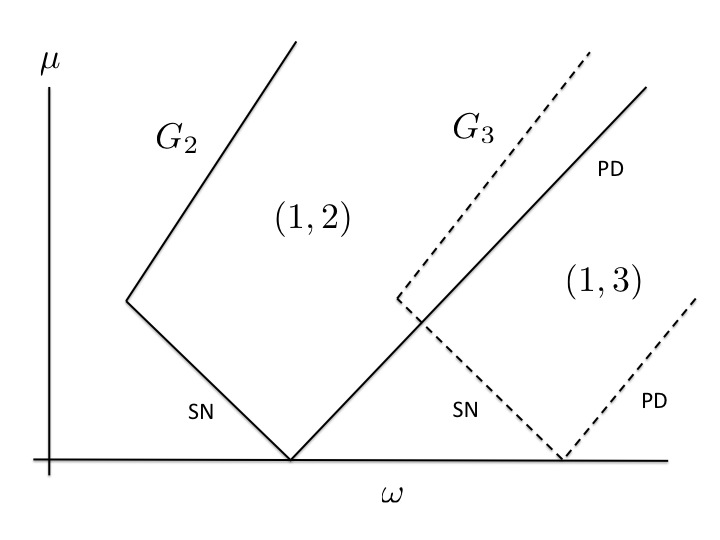}
\caption{A schematic of the above figure showing the domains of existence of the $(1,2)$ and the $(1,3)$ orbits. In this figure we see the grazing curves $G_n$ and the curves of saddle-node and period-doubling bifurcations. }
\label{leafschem}
\end{figure}

\begin{figure}[htb!]
\centering
\includegraphics[scale=0.18]{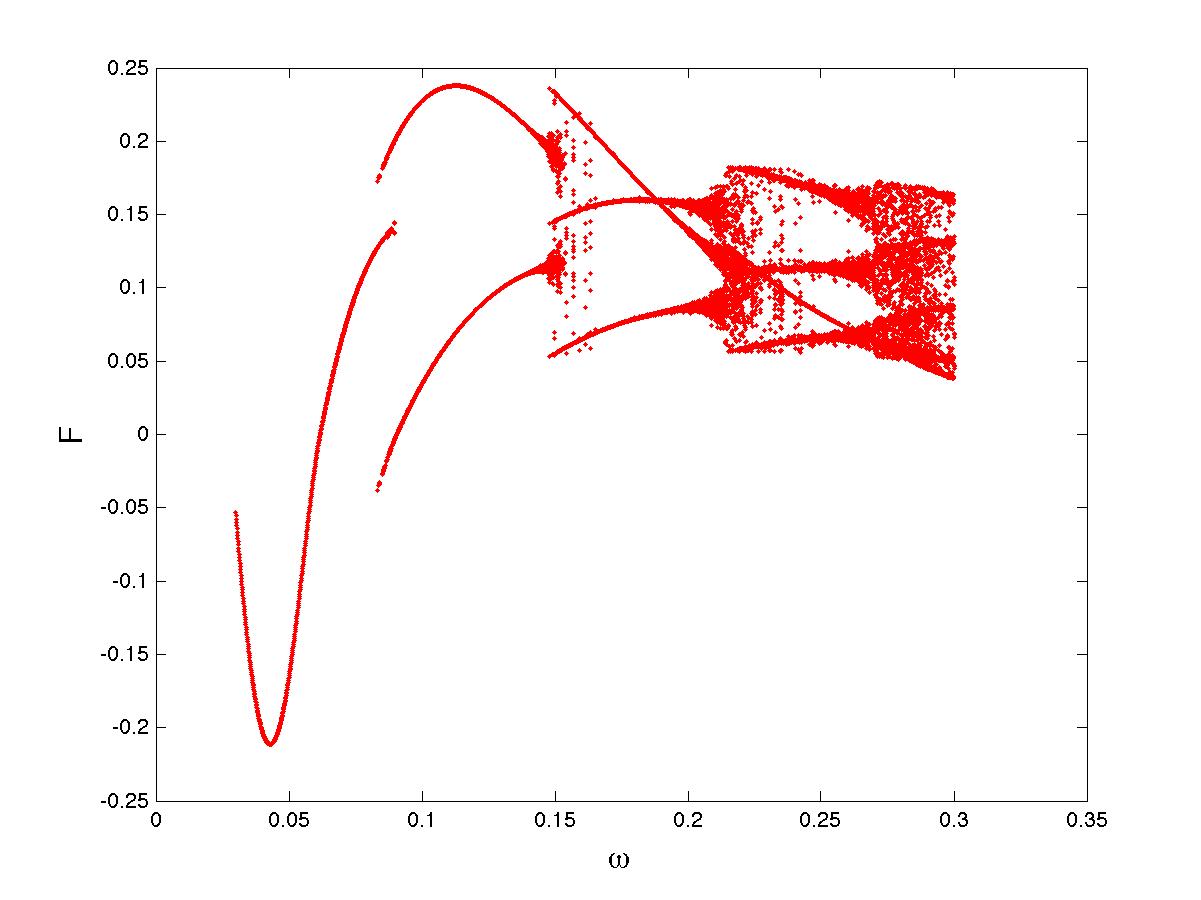}
\caption{A Mont\'e-Carlo bifurcation diagram when $\mu = 1$ showing the periodic and other orbits, and the grazing and other bifurcations. }
\label{bibbif}
\end{figure}

\noindent We calculate the curves $G_n(\omega)$ by fixing $\mu$ and following the $(1,n)$ orbits as $\omega$ decreases, stopping at $G_n(\mu)$ when we get the additional grazing impact. The curve $G_3$ is
also plotted (in solid mauve) on  Figure \ref{fig:leaf} together with the curve of saddle-node bifurcations of the $(1,3)$ orbit. We can see clearly 
from this figure that the tongue for the existence of each 
$(1,n)$ orbit is bounded on the {\em left} by the curves $G_n$. In particular if $\omega$ decreases through the value $G_n(\mu)$ we see a {\em sudden change} from a $(1,n)$ orbit to a $(1,n-1)$ orbit. In Figure \ref{leafschem} we give a schematic of the regions of existence. In this schematic we see the lines of saddle-node and period-doubling bifurcations emerging from the resonance points when $\mu=0$ and $\omega = n \omega^*$. The curve of saddle-node bifurcations is then 'destroyed' when it intersects the lines $G_n$ of grazing bifurcations. These lines then become the leftmost boundary of the domains of existence. Finally in Figure \ref{bibbif} we give a bifurcation diagram computed using the Mont\'e-Carlo method which fixes $\mu = 1$ and takes
$0.03 < \omega < 0.3$. In this diagram we can see the $(1,1), (1,2), (1,3)$ and $(1,4)$ periodic orbits, together with the grazing bifurcations at $\omega_g = 0.084, 0.148, 0.21$, the period-doubling bifurcations at $\omega_{PD} = 0.09, 0.16, 0.22$ and regions of chaotic and quasi-periodic behaviour.

\subsection{Location of the grazing lines}

\noindent We now make some estimates of the location of the curves $G_n$. We firstly show that $G_n$ does not exist for small values of $\mu$.

\vspace{0.1in}

\noindent {\bf Lemma 8.2} {\em If $r ^{\pm} \equiv -c^T L^{-1} b^{\pm} + d$ is bounded away from zero, and if $\mu \ll |r^{\pm}|$ then we cannot have a grazing impact.}

\vspace{0.1in}

\noindent {\bf Proof} We firstly note that for very small $\mu$ and for $\omega$ close to $n \omega^*$ the observed periodic orbits will be ${\cal O}(\mu)$  perturbations of the 'unforced' periodic orbit at $\mu = 0$ of period $2\pi/\omega^*$. This orbit is shown in \cite{morupisi2020analysis} to arise at a border collision bifurcation. The unforced orbit  is made up of piece-wise decaying exponential curves on which $F$ converges to virtual limits given by
\begin{equation}
    r ^{\pm} = -c^T L^{-1} b^{\pm} + d
    \label{res1}
\end{equation}
so that there are non-zero constants $A^{\pm}$ so that
$$F(t) \approx A^{\pm} e^{-\lambda_1 t} + r^{\pm}, \quad \dot{F}(t) \approx -\lambda_1 A^{\pm} e^{-\lambda_1 t},$$
where $\lambda_{min}$ is the smallest (in magnitude) eigenvalue of $L$.
It follows that if $F(t) = 0$ then $\dot{F}(t) \approx \lambda_1 r^{\pm}$ and there is no possible graze of the unforced periodic solution. 
If the unforced system is now perturbed by a periodic function of amplitude proportional to $\mu$, then as we have seen the perturbed problem will have a periodic solution which differs from the unforced periodic function by a perturbation also of ${\cal O}(\mu)$.
As the unforced periodic solution does not have a grazing impact, then the perturbed system will have a nearby impact, for which 
$$\dot{F} \approx \lambda_1 r^{\pm} + {\cal O}(\mu).$$
It follows that, provided that $|r^{\pm}| > 0$ and  $\mu \ll |r^{\pm}|$,
then this cannot be a grazing impact. \qed

\vspace{0.1in}

\noindent We note from this lemma, that if $\mu$ is small then all intersections are transversal, and the dynamics of the PP04 model is completely described by that of smooth dynamical systems. This explains the 'regular form' Figure \ref{fig:leaf} for small values of $\mu.$ See \cite{morupisi2020analysis} for an exploration of this dynamics. 

\vspace{0.1in}

\noindent However, if we fix $\omega = \omega_g$ then as $\mu$ increases there as the periodic in time perturbation to the unforced periodic orbit is proportional to $\mu$, this will eventually dominate the unforced period solution. Hence, as $\mu$ increases there will be a {\bf first} $\mu_g > 0$ at which a grazing impact occurs for which we have $\omega_g = G_n(\mu_g)$.

\vspace{0.1in}

\noindent We note further from the theory presented earlier that
the perturbation is of the form
$$ {\mathbf a} \cos(\omega t) + {\mathbf b} \sin (\omega t)$$
where 
$${\mathbf a} = \mu \omega (L^2 + \omega^2)^{-1} {\mathbf e}, \quad {\mathbf b}
= \mu L (L^2 + \omega^2)^{-1}{\mathbf e} .$$
For large $\omega/\lambda_1$ the magnitude $m$ of the perturbation is proportional to $\mu/\omega$. Hence we expect to see lines of grazing bifurcations in which $\omega_g$ is approximately proportional to $\mu_g$ when both are large. This is what is observed in Figure \ref{fig:leaf}.

\vspace{0.1in}

\noindent Whilst we have plotted the regions of existence in the $(\omega, \mu)$ plane, the arguments that we have presented apply to perturbations of any of the other parameters in the PP04 model. Indeed, as we saw from the numerical calculations presented in Section 3, we can see very similar transitions when the parameter $d$ is varied, with regions of existence of the periodic orbits again bounded by lines of grazing, and period-doubling, bifurcations.

\section{Quasi-periodic perturbations}
\noindent
 
 \noindent The earlier analysis has concentrated on periodic insolation forcing. However, in a true climate situation, the insolation forcing (due to the Milankovich cycles) is quasi-periodic.  This will introduce a  perturbation to the periodic orbits that we have considered in the previous sections. In \cite{mitsui2014dynamics}, \cite{mitsui2015bifurcations} the impact on PP04 when considered as a smooth dynamical system was discussed and the possible existence of strange non-chaotic attractors was described. Mitsui noted, however, that it would be more appropriate to consider interplay between the non-smooth structure of the model and the quasi-periodic perturbation.  We now consider this case briefly, and a fuller discussion will be given in a subsequent paper.  
 
 \vspace{0.1in}
 
 \noindent If a periodic orbit of the periodically forced system  is not close to grazing, with all impacts transversal, the dynamical system described by the PP04 model close to this orbit is locally smooth. Hence any small perturbations of the periodic orbit will be regular. See, for example, the calculations reported in \cite{ashwin2015middle} and \cite{morupisi2020analysis}. In contrast, we expect that the perturbation of a near grazing periodic orbit by the quasi-periodic terms will be much more pronounced due to the effect of the discontinuities in the associated Poincar\'e map. 
 
 \vspace{0.1in}
 
\begin{figure}[htb!]
\centering
\includegraphics[scale=0.3]{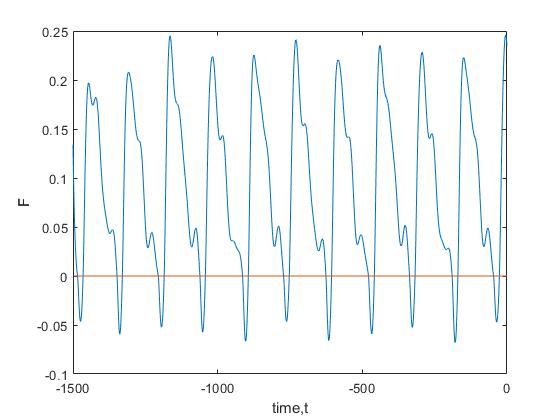}
\includegraphics[scale=0.3]{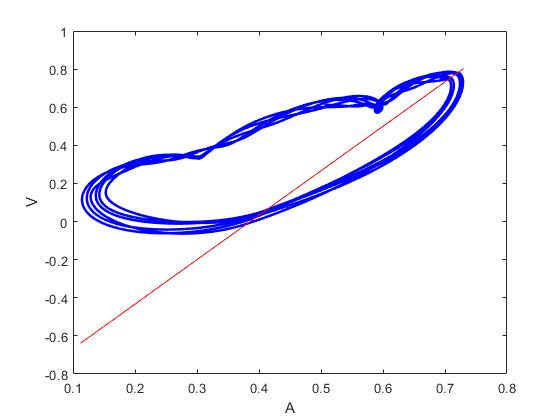}
\caption{(left) The function $F(t)$ obtained  when $\mu_{1} = 0.3,\omega_{1}=0.13$ and $\mu_{2}=0.2,\omega_{2}= 0.146$ showing the quasi-periodic orbit when a period $(1,3)$ periodic orbit which is not close to grazing is perturbed. We see that the quasi-periodic orbit lies close to the periodic orbit. (right) the same orbit in the $A-V$ phase plane}
\label{peridcpert}
\end{figure}

\begin{figure}[htb!]
\centering
\includegraphics[scale=0.3]{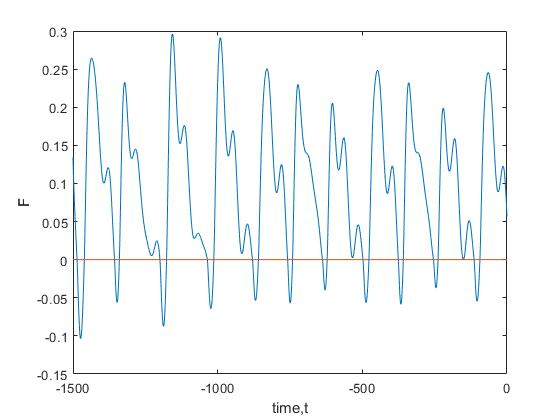}
\includegraphics[scale=0.3]{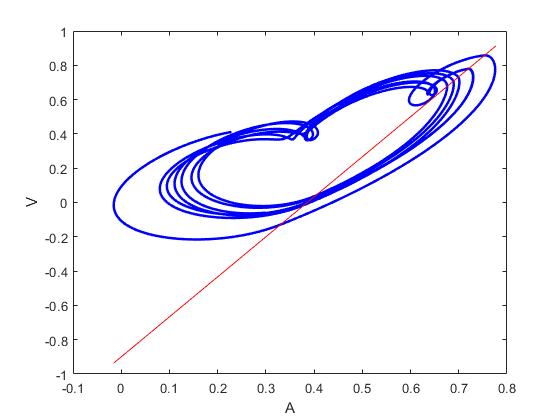}
\caption{(left) The function $F(t)$ obtained  when $\mu_{1} = 0.3,\omega_{1}=0.115$ and $\mu_{2}=0.25,\omega_{2}= 0.1476$ 
(right) the same orbit in the $A-V$ phase plane. In this case we see that the quasi-periodic orbits differs significantly from the periodic orbit.}

%(right) $\mu_{1} = 0.3,\omega_{1}=0.115$ and %$\mu_{2}=0.1,\omega_{2}= 0.126$, showing that the quasi-periodic %orbit differs significantly from the periodic orbit.}
\label{pertng1}
\end{figure}

 \noindent As an example we consider the case of $\mu = 0.3$ and take (a) $\omega  = 0.13$ and (b) $\omega = 0.115$. From the earlier calculations we know that the first case has an orbit which is not close to grazing, and in the second case the orbit is close to grazing.  We perturb the  periodic orbit by adding another forcing frequency to the original. To do this we set
the insolation forcing to equal the quasi-periodic function:
$$f(t) = \mu_1 \sin(\omega_1 t) + \mu_2 \sin(\omega_2 t)$$
The resulting dynamics when we take $\omega_1 = 0.13, \mu_1 = 0.3, \omega_2 = 0.1476, \mu_2 = 0.25$ is shown in Figure  \ref{peridcpert} where we show both the time evolution of $F(t)$ and the dynamics in the $A-V$ phase-plane.  As expected that the perturbed quasi-periodic orbit is close to the original periodic orbit. In contrast, in Figure \ref{pertng1} we take $\omega_1 = 0.115, \mu_1 = 0.3, \omega_2 = 0.1476, \mu_2 = 0.25$. The quasi-periodic orbit now deviates significantly from the periodic orbit, with transitions observed at the near grazing events. The latter Figure has the general structure of the time-series observed from the reconstructed paleo-climate data.

 \noindent A fuller study of the dynamics which arises when the PP04 model is forced quasi-periodically (including a fuller discussion of SNAs) will be the subject of a forthcoming paper.

\section{Conclusions}
\noindent In this paper we have studied the type of dynamics and  transitions observed between  different types of periodic orbits for the periodic forced  PP04 model. This has  revealed  the existence of regions of synchronised periodic solutions $(1,n)$ orbits. For smaller values of $(\mu,\omega)$ the boundaries of the existence regions are given by saddle-node or period doubling bifurcations. However, for larger values of $\mu$ the synchronised periodic solutions can lose stability at grazing events. Such grazing bifurcations are highly destabilising and lead to sudden, and dramatic, changes in behaviour. These can be a sudden shift in the type of periodic orbit seen, or a sudden change in the phase of the orbit. The sets of the grazing bifurcations then mark the boundaries both of the regions of existence of the periodic solutions as parameters are varied, and also of the domains of attraction of these orbits. The resulting transitions between the periodic orbits when either external parameters (such as $\omega$) or internal parameters (such as $d$) are varied,  resembles the changes that  we associate with the MPT with rapid transitions between different periodic states, and in particular a change from a low amplitude $(1,2)$ orbit to a higher amplitude $(1,3)$ orbit. We suggest from this evidence, and the generic nature of grazing,  that the MPT may have arisen from a grazing bifurcation under a parameter change. 

\vspace{0.1in}

 \noindent The  astronomical forcing  is quasi-periodic and therefore we briefly considered the impact of a quasi-periodic forcing on the PP04 model. We have shown that grazing events can cause significant changes to the periodic orbits when the extra quasi-periodic forcing is included.  More work is needed to understand very rich dynamics of the quasi-periodically forced PP04 model, and this will form the basis of a further study.

\vspace{0.2in}

\noindent {\bf Acknowledgements} We would like to thank the two anonymous referees for their careful reading, and comments on an earlier version of this paper, which has led to a significant improvement in the results presented. This research was funded in part by an award from the Botswana International University of Science and Technology (BIUST).
\bibliographystyle{plain}

%bibliography{references}
\end{document}